\documentclass[12pt]{elsarticle}
\usepackage{amsthm}
\usepackage{amsmath}
\usepackage{here}
\usepackage{multirow}
\usepackage[svgnames]{xcolor}
\usepackage{comment}
\usepackage{enumerate}
\usepackage[top=20truemm,bottom=20truemm,left=20truemm,right=20truemm]{geometry}
\usepackage{ulem}

\usepackage{bm}
\theoremstyle{remark}
\newtheorem{remark}{Remark}

\begin{document}

\begin{frontmatter}

\title{Symmetry-preserving enforcement of low-dissipation method based on boundary variation diminishing principle}

\author[a]{Hiro Wakimura\corref{cor1}}
\ead{wakimura.h.aa@m.titech.ac.jp, ijsp090724@gmail.com}
\cortext[cor1]{Corresponding author}

\author[a]{Shinichi Takagi}
\author[a]{Feng Xiao\corref{cor1}}
\ead{xiao.f.aa@m.titech.ac.jp}

\address[a]{Department of Mechanical Engineering, Tokyo Institute of Technology, 2-12-1 Ookayama, Meguro-ku, Tokyo, 152-8550, Japan}

\begin{abstract}
A class of high-order shock-capturing schemes, P$_n$T$_m$-BVD (Deng et al., J. Comp. Phys., 386:323-349, 2019; Comput. \& Fluids, 200:104433, 2020.) schemes, have been devised to solve the Euler equations with substantially reduced numerical dissipation, which enable high-resolution simulations to resolve flow structures of wider range scales. In such simulations with low dissipation, errors of round-off level might grow and contaminate the numerical solutions. A typical example of such problems is the loss of symmetry in the numerical solutions for physical problems of symmetric configurations even if the schemes are mathematically in line with the symmetry rules. In this study, the mechanisms of symmetry-breaking in a finite volume framework with the P$_4$T$_2$-BVD reconstruction scheme are thoroughly examined. Particular attention has been paid to remove the possible causes due to the lack of associativity in floating-point arithmetic which is associated with round-off errors. 
Modifications and new techniques are proposed to completely remove the possible causes for symmetry breaking in different components of the P$_4$T$_2$-BVD finite volume solver. Benchmark tests that have symmetric solution structures are used to verify the proposed methods. The numerical results demonstrate the perfect symmetric solution structures.

\end{abstract}


\end{frontmatter}


\section{Introduction}

The most pronounced feature of compressible flow is that both smooth and discontinuous flow structures, such as vortices, acoustic waves, turbulence, shock waves, and contact surfaces, evolve simultaneously in the computational domain of interest. It remains a big challenge to design numerical schemes that capture both discontinuous and smooth solutions of wide-range scales. The finite volume method based on the Godunov scheme \cite{Godunov1959} is the most widely-used mainstream numerical framework for compressible flows because of its rigorous conservativeness at the discrete level, which is crucial to correctly capture the nonlinear discontinuous solutions. Tremendous efforts have been made to develop high-resolution schemes for spatial reconstruction in the finite volume framework during the past decades. In general, spatial reconstructions using high-order polynomials are needed to reduce numerical dissipation so as to resolve wide-range scale flow structures, however, as stated in Godunov’s theorem \cite{Godunov1959}, any linear scheme with a polynomial equal to or higher than 1st-order tends to generate numerical oscillations in the vicinity of the discontinuities. Hence, a large class of reconstruction schemes have been extensively explored based on the so-called nonlinear limiting projection. An ideal high-resolution shock-capturing scheme is expected to be able to effectively suppress spurious oscillations while retaining high order and reducing excessive numerical dissipation or viscosity to adequately resolve flow structures of small scales. 

TVD (Total Variation Diminishing) schemes \cite{Harten1983} including MUSCL (Monotone Upstream-centered Schemes for Conservation Law) scheme \cite{VanLeer1977} maintains monotonicity by introducing the slope limiter, which suppresses the numerical oscillations near the discontinuities while retrieving nearly 2nd-order accuracy for smooth solutions. Although the numerical dissipation of the TVD schemes is rather smaller than that of the 1st-order Godunov scheme, excessive numerical dissipation is still too large to make them attractive in applications where vortices and turbulence are of great importance.  Progress has been made to use high-order polynomials with nonlinear limiting projection. Being the most representative schemes, ENO (Essentially Non-Oscillatory) scheme \cite{Harten1987,Shu1988} and WENO (Weighted Essentially Non-Oscillatory) scheme \cite{Liu1994,Jiang1996} were proposed. These schemes achieved high convergence order (higher than 2nd-order) while suppressing numerical oscillation by adaptively selecting or weighting the smooth sub-stencils. Successive researches were carried out to devise variants of the smooth indicator   \cite{Henrick2005,Borges2008,Ha2013,Fan2014,Kim2016,Fu2016,Acker2016,Fu2017,Rathan2018,Fu2018} to suppress the numerical dissipation near the discontinuities. Nevertheless, the numerical dissipation of these schemes remains non-negligible and tends to smear out the fluid structures such as contact discontinuities and small-scale vortices. Moreover, the nonlinear limiting projection, which degrades the degree of reconstruction function near a discontinuity, can hardly recover the low-dissipation property of the original linear scheme using the unlimited polynomial with constant coefficients.

Based on the observation that the dissipation term in the approximate Riemann solvers can be effectively reduced by minimizing the difference between the reconstructed left- and right-side values at cell boundaries (referred to as BV (boundary variation)), a new class of reconstruction schemes so-called BVD (Boundary Variation Diminishing) schemes have been proposed in \cite{Sun2016,Deng2017,Xie2017a,Deng2018a,Deng2018,Deng2019,Tann2019,Deng2020,Tann2020,Cheng2021,Jiang2021}. The first BVD scheme \cite{Sun2016}, so-called WENO-THINC-BVD scheme, switches the reconstruction function to WENO scheme for smooth solutions and THINC (Tangent of Hyperbola for INterface Capturing) for discontinuous solutions. The THINC scheme was originally developed for capturing moving interfaces \cite{Xiao2005,Xiao2011,Xie2017} and thus more suitable for representing the discontinuous solutions of a step-jump shape. As the THINC can more accurately mimic a discontinuous solution within a mesh cell, it effectively reduces the BV of the reconstructed values across cell boundaries, which enables the BVD criterion to choose the THINC function in presence of a discontinuity. Meanwhile, the polynomial-based WENO interpolation leads to smaller BV values for smooth solutions, so the BVD scheme will select the WENO scheme as the reconstruction function. As a result, the WENO-THINC-BVD scheme significantly reduces the numerical dissipation and shows superiority in resolving both smooth and discontinuous solutions. As a high-resolution shock-capturing scheme without using the conventional nonlinear limiting projection, the P$_4$T$_2$-BVD (polynomial of 4-degree and THINC function of 2-level reconstruction based on BVD algorithm) scheme was proposed in \cite{Deng2019}, which combines the unlimited polynomial of the four degree and THINC schemes, THINC($\beta_s$) and THINC($\beta_l$) with different steepness parameters.  The P$_4$T$_2$-BVD scheme uses a 2-stage BVD algorithm to select the reconstruction function under the BVD principle, and greatly improves numerical solutions regarding numerical dissipation. The numerical results of the P$_4$T$_2$-BVD scheme can resolve both smooth and discontinuous solutions with superior fidelity compared to other existing shock-capturing schemes that use the nonlinear limiters to the polynomial interpolants in the presence of discontinuities. It is also demonstrated in \cite{Deng2020} that higher-order BVD schemes, P$_n$T$_m$-BVD  (polynomial of $n$-degree and THINC function of $m$-level reconstruction based on BVD algorithm) schemes, can be designed by using high-order unlimited polynomials and THINC functions with multi-stage BVD algorithms. It is noted that the BVD principle is general and can be applied to design low dissipation schemes using other candidate functions for reconstruction \cite{Chamarthi2021} or other numerical framework \cite{Ruan2020}.  With substantially reduced numerical dissipation, the BVD schemes can effectively preserve small-scale flow structures from being dampened out.

As discussed above, shock-capturing schemes evolve to meet the requirement from applications where numerical dissipation has to be substantially reduced. 
In practice, numerical dissipation can be reduced by either designing numerical methods of less dissipation error or refining computational grids for simulations. Modern computer hardwares facilitate the use of advanced high-order schemes with fine mesh resolution, which enables us to produce simulation results with extremely low numerical dissipation.  
With less numerical dissipation, flow structures of small scale remain unattenuated, so do perturbations stemming from floating-point arithmetic errors that are unavoidable in numerical processing. As a result, the numerical solution may lose spatial symmetry even started from a symmetrical setup. It is commonly observed that high-resolution simulations with fine meshes usually lead to undesirable asymmetrical numerical results \cite{Shi2003,Remacle2003,Liska2003,Ha2005,Fu2016,Zhao2018,Zhang2020,Li2020,Ruan2020,Peng2021,Li2021,Tann2020} in application problems that possess physical symmetry in solution structures, such as the Rayleigh-Taylor instability (RTI) and the symmetric implosion benchmark problems. The symmetry-breaking results might not be visible in simulations with significant numerical dissipation, and have not been seen as a serious problem in the past.

Remacle et al. \cite{Remacle2003} initially reported that the symmetry error is caused by the rounding error of floating-point arithmetic, which grows asymmetrically with time evolution. Don et al. suggested a numerically stable form of the smooth indicator in the WENO framework \cite{Don2018,Don2020}. In this scheme, the symmetry error was effectively reduced compared to the original 7th- and 9th-order WENO schemes, but the cause of the symmetry error has not been completely eliminated.  Wang et al. \cite{Wang2020} preserved the symmetry property by replacing the two values of solutions at symmetrically placed cells with their average value. Fleischmann et al. \cite{Fleischmann2019} attribute the cause of asymmetric rounding errors to the lack of associativity in summations or multiplications for more than two components and carefully analyzed the solution procedure of some WENO/TENO schemes. It is shown that by adjusting the orders of summations and multiplications throughout the calculation procedure, the symmetry errors can be completely eliminated in the high-resolution simulations that implement WENO-like schemes on fine meshes. To our knowledge, \cite{Fleischmann2019} is the first work that implements  the WENO-type schemes to ensure exact symmetry property in numerical solutions.

In this paper, we propose symmetry-preserving P$_4$T$_2$-BVD scheme which belongs to another sort of high-order reconstruction approach different from the WENO-type schemes dealt with in \cite{Fleischmann2019}. As the original THINC scheme may generate the symmetry error, a new formulation of symmetry-preserving THINC scheme is presented. Moreover, we thoroughly examined the possible causes in the numerical formulations for different components of the numerical solver, including the transformation between conservative and characteristic variables, spatial reconstruction and approximate Riemann solver.  We clarified some important issues which have not been specified in \cite{Fleischmann2019} or other existing literature. Countermeasures are proposed to exactly ensure the spatial symmetry in numerical solutions to physical benchmarks of symmetrical configurations. Not limited to the BVD method,  the symmetry-preserving techniques introduced in this paper are efficient and easy to be embedded into the  finite volume framework. Essentially, the present symmetry-preserving techniques modify the operations that might generate different machine errors without compromising the performances of the original schemes.

The rest of this paper are organized as follows. In section \ref{sec:Governing equation}, the inviscid Euler equations are briefly reviewed as the governing equations. In section \ref{sec:Numerical methods}, the calculation procedure of P$_4$T$_2$-BVD scheme is briefly explained. In section \ref{sec:Symmetry}, we thoroughly examine the mechanisms of the symmetry-error generation  in the finite volume framework using  P$_4$T$_2$-BVD reconstruction and propose symmetry-preserving techniques to eliminate the possible causes for symmetry-breaking. Numerical results of benchmark tests which have symmetry properties are shown in section \ref{sec:Results} for verification.  Conclusion is given in section \ref{sec:Conclusion} to end this paper.

\section{Governing equations}
\label{sec:Governing equation}
We consider inviscid compressible flow, and use the following 2D Euler equations as the governing equations.
\begin{equation}
\frac{\partial \mathbf{U}}{\partial t}+\frac{\partial \mathbf{F(U)}}{\partial x}+\frac{\partial \mathbf{G(U)}}{\partial y}=0,
\label{euler}
\end{equation}
where $\mathbf{U}$ is the vector of conservative variable, and $\mathbf{F}$ and $\mathbf{G}$ are the vectors of flux functions in $x$- and $y$-directions respectively. The components of these vectors are given explicitly as below,
\begin{equation}
\mathbf{U}=
\begin{pmatrix}
\rho \\ \rho u \\ \rho v \\ E \\
\end{pmatrix}
,\ \mathbf{F(U)}=
\begin{pmatrix}
\rho u \\ \rho u^2+p \\ \rho uv \\ (E+p)u \\
\end{pmatrix}
,\ \mathbf{G(U)}=
\begin{pmatrix}
\rho v \\ \rho vu \\ \rho v^2+p \\ (E+p)v \\
\end{pmatrix}
,
\label{flux_function}
\end{equation}
where $\rho$ is the density, $u$ and $v$ are the $x$- and $y$-components of the velocity, $E$ is the total energy per unit volume, and $p$ the static pressure.

The Euler equations \eqref{euler} are closed by adding the equation of state (EOS) of the ideal gas as follows,
\begin{equation}
p=(\gamma-1)\left( E-\frac{1}{2}\rho \left( u^2+v^2\right)\right),
\end{equation}
where $\gamma$ is the specific heat ratio, and $\gamma=1.4$ is used in the present work unless specifically noted.

In the numerical formulation presented in this work, the spatial reconstructions are implemented in terms of characteristic variables, and the arrangement of arithmetic operations of transformation among conservative variables and characteristic variables is crucial to maintain the symmetry in numerical solutions. Thus, we provide explicitly the details of the transformation relations,  
\begin{align}
&\mathbf{W}_x=\mathbf{L}_x\cdot \mathbf{U} \ \ \text{and} \ \ \mathbf{U}=\mathbf{R}_x\cdot \mathbf{W}_x, \label{trans-x}\\
&\mathbf{W}_y=\mathbf{L}_y\cdot \mathbf{U} \ \ \text{and} \ \ \mathbf{U}=\mathbf{R}_y\cdot \mathbf{W}_y, \label{trans-y}
\end{align}
where $\mathbf{W}_x$ and $\mathbf{W}_y$ are characteristic variables used for reconstructions in $x$- and $y$-directions, and denoted respectively by  
\begin{align}
\mathbf{W}_x=
\begin{pmatrix}
w^{(u-c)} \\ w^{(u)} \\ w^{(u+c)} \\ w^{(u\perp)}
\end{pmatrix}
,\ \mathbf{W}_y=
\begin{pmatrix}
w^{(v-c)} \\ w^{(v)} \\ w^{(v+c)} \\ w^{(v\perp)}
\end{pmatrix}.
\end{align}
Here, $c=\sqrt{\frac{\gamma p}{\rho}}$ is the sound speed, and symbol  ``$\perp$'' is used to distinguish the components corresponding to the two eigenvalues $u$ and $v$. The left eigen matrices and the right eigen matrices in $x$- and $y$-directions are written as
\begin{align}
&
\begin{aligned}
\mathbf{L}_x=
\begin{pmatrix}
\vec{l}^{(u-c)} \\ \vec{l}^{(u)} \\ \vec{l}^{(u+c)} \\ \vec{l}^{(u\perp)}
\end{pmatrix}
=
\begin{pmatrix}
\frac{1}{2}\left(b_1+\frac{u}{c}\right) & -\frac{1}{2}\left(\frac{1}{c}+b_2u\right) & -\frac{1}{2}b_2v & \frac{1}{2}b_2 \\
1-b_1 & b_2u & b_2v & -b_2 \\
\frac{1}{2}\left(b_1-\frac{u}{c}\right) & \frac{1}{2}\left(\frac{1}{c}-b_2u\right) & -\frac{1}{2}b_2v & \frac{1}{2}b_2 \\
-v & 0 & 1 & 0
\end{pmatrix}=
\begin{pmatrix}
\vec{l}_{x1} & \vec{l}_{x2}&\vec{l}_{x3} & \vec{l}_{x4} 
\end{pmatrix}
,\\
\mathbf{L}_y=
\begin{pmatrix}
\vec{l}^{(v-c)} \\ \vec{l}^{(v)} \\ \vec{l}^{(v+c)} \\ \vec{l}^{(v\perp)}
\end{pmatrix}
=
\begin{pmatrix}
\frac{1}{2}\left(b_1+\frac{v}{c}\right) & -\frac{1}{2}b_2u & -\frac{1}{2}\left(\frac{1}{c}+b_2v\right) & \frac{1}{2}b_2 \\
1-b_1 & b_2u & b_2v & -b_2 \\
\frac{1}{2}\left(b_1-\frac{v}{c}\right) & -\frac{1}{2}b_2u & \frac{1}{2}\left(\frac{1}{c}-b_2v\right) & \frac{1}{2}b_2 \\
-u & 1 & 0 & 0
\end{pmatrix}=\begin{pmatrix}
\vec{l}_{y1} & \vec{l}_{y2}&\vec{l}_{y3} & \vec{l}_{y4} 
\end{pmatrix}
,
\label{eigen_L}
\end{aligned}
\\
&
\begin{aligned}
\mathbf{R}_x=
\begin{pmatrix}
\vec{r}^{(u-c)} & \vec{r}^{(u)} & \vec{r}^{(u+c)} & \vec{r}^{(u\perp)}
\end{pmatrix}
=
\begin{pmatrix}
1 & 1 & 1 & 0 \\
u-c & u & u+c & 0 \\
v & v & v & 1 \\
H-uc & \frac{u^2+v^2}{2} & H+uc & v
\end{pmatrix}=\begin{pmatrix}
\vec{r}_{x1} & \vec{r}_{x2}&\vec{r}_{x3} & \vec{r}_{x4} 
\end{pmatrix}
,\\
\mathbf{R}_y=
\begin{pmatrix}
\vec{r}^{(v-c)} & \vec{r}^{(v)} & \vec{r}^{(v+c)} & \vec{r}^{(v\perp)}
\end{pmatrix}
=
\begin{pmatrix}
1 & 1 & 1 & 0 \\
u & u & u & 1 \\
v-c & v & v+c & 0 \\
H-vc & \frac{u^2+v^2}{2} & H+vc & u
\end{pmatrix}=\begin{pmatrix}
\vec{r}_{y1} & \vec{r}_{y2}&\vec{r}_{y3} & \vec{r}_{y4} 
\end{pmatrix}
,
\label{eigen_R}
\end{aligned}
\end{align}
where $b_1=\frac{u^2+v^2}{2}\frac{\gamma-1}{c^2},\ b_2=\frac{\gamma-1}{c^2}$, and the enthalpy $H=\frac{E+p}{\rho}$. It is noted that the eigen vectors in \eqref{eigen_L} and \eqref{eigen_R} are written in the natural order. However, the arrangement of eigenvectors in the eigen matrices affects the symmetry property in numerical solution, which will be discussed later. 

\section{Numerical methods}
\label{sec:Numerical methods}
For the sake of simplicity, we use the 1D scalar hyperbolic conservation law as the model equation to describe numerical schemes, which is cast in the general form as,
\begin{equation}
\frac{\partial q}{\partial t}+\frac{\partial f(q)}{\partial x}=0,
\label{scalar_cons}
\end{equation}
where $q(x,t)$ is the solution of the conservative variable and $f(q)$ is the flux function. The characteristic speed $a=\partial f(q)/\partial q$ is a real number because of the hyperbolicity of Eq. \eqref{scalar_cons}.

\subsection{Finite volume method}

We divide the computational domain into $N$ non-overlapping cells, $\Omega_i=[x_{i-1/2}, x_{i+1/2}],i=1,2,...,N$, and assume that the cell size $\Delta x=x_{i+1/2}-x_{i-1/2}$ is uniform over the computational domain for brevity.

For a standard finite volume method, the volume-integrated average of the solution in cell $\Omega _i$ is defined as,
\begin{equation}
\bar{q}_i(t)\equiv \frac{1}{\Delta x}\int_{x_{i-1/2}}^{x_{i+1/2}}q(x,t)dx.
\label{VIA}
\end{equation}
The solution $\bar{q}_i(t)$ for each cell $\Omega_i$ is updated by the semi-discrete version of Eq.  \eqref{scalar_cons} as
\begin{equation}
\frac{d \bar{q}_i}{d t}=-\frac{1}{\Delta x}\left( \hat{f}_{i+\frac{1}{2}}-\hat{f}_{i-\frac{1}{2}}\right),
\label{semi-discrete}
\end{equation}
where $\hat{f}_{i+\frac{1}{2}}$ is the numerical flux at cell boundary $x_{i+\frac{1}{2}}$. In a Godunov finite volume method, numerical flux  $\hat{f}_{i+\frac{1}{2}}$ is calculated by an approximate Riemann solver,
\begin{equation}
\hat{f}_{i+\frac{1}{2}}=f^{\rm Riemann}\left( q_{i+\frac{1}{2}}^L, q_{i+\frac{1}{2}}^R\right),
\label{Riemann_solver}
\end{equation}
where $q_{i+1/2}^L$ and $q_{i+1/2}^R$ stand for the left-side and right-side values of solution $q$ at the cell boundary $x_{i+1/2}$ computed from spatial reconstructions.
Although various  approximate Riemann solvers have been  proposed in different forms, we present them in a canonical form that consists of a central scheme and a dissipation (or viscosity) term, 
\begin{equation}
f^{\rm Riemann}\left( q_{i+\frac{1}{2}}^L, q_{i+\frac{1}{2}}^R\right)=\frac{1}{2}\left( f\left( q_{i+\frac{1}{2}}^L\right)+f\left( q_{i+\frac{1}{2}}^R\right)\right)-\frac{1}{2}\left| a_{i+\frac{1}{2}}\right| \left( q_{i+\frac{1}{2}}^R-q_{i+\frac{1}{2}}^L\right).
\label{Riemann_solver_canonical}
\end{equation}
Given a Riemann solver, the spatial reconstruction scheme that calculates the boundary values $q_{i+1/2}^L$ and $q_{i+1/2}^R$ over the left- and right-biased stencils plays the key role to determine the solution quality of numerical results. In next subsection, we focus on the numerical scheme, so-called the P$_4$T$_2$-BVD scheme, for spatial reconstruction.

\subsection{The P$_4$T$_2$-BVD scheme}

BVD principle \cite{Sun2016} has been proposed as a general guideline for  spatial reconstruction with reduced numerical dissipation.
In the canonical form of the Riemann flux shown in Eq. \eqref{Riemann_solver_canonical}, the second term on the right-hand side can be interpreted as numerical viscosity/dissipation. It is observed that the difference of the reconstructed  values at cell boundary, so-called boundary variation (BV), 
\begin{equation}
BV_{i+\frac{1}{2}}=\left| q_{i+\frac{1}{2}}^R-q_{i+\frac{1}{2}}^L\right|, 
\end{equation}
 is in proportion to numerical dissipation. The BVD principle suggests that spatial reconstruction should be devised so that the value of $BV_{i+1/2}$ is minimized. In BVD schemes, the BVD principle is used as a criterion to select a reconstruction function from multiple candidate functions. Some practical BVD selection algorithms have been proposed so far to design variants of BVD schemes \cite{Sun2016,Deng2018a,Xie2017a,Deng2018}. High-fidelity schemes which have superior properties in suppressing both numerical oscillation and dissipation can be designed by using  BVD algorithms and properly chosen candidate interpolation functions.

Being a representative BVD scheme, the P$_4$T$_2$-BVD scheme \cite{Deng2019,Deng2020}, makes use of the unlimited 4th-degree polynomial and two THINC functions of different steepness parameters, namely THINC($\beta_s$) and THINC($\beta_l$) . The selection of the reconstruction function is performed by comparing the values of the total $BV$ ($TBV$) \cite{Deng2018} of the reconstructed values in a 2-stage algorithm. It is demonstrated that the P$_4$T$_2$-BVD scheme can effectively remove numerical oscillations without using the conventional nonlinear limiting projections which are widely used in other high-resolution shock capturing schemes. For completeness, we describe  the P$_4$T$_2$-BVD scheme in this subsection.

\subsubsection{Candidate interpolant 1: 4th-degree polynomial scheme}

In the P$_4$T$_2$-BVD scheme, the physical quantity $q$ is reconstructed by the linear 4th-degree polynomial as one of the candidate interpolation functions. The reconstruction function $\mathcal{Q}(x)$ for target cell $\Omega_i$ is written as
\begin{equation}
\mathcal{Q}_i^{P_4}(x)=\sum_{k=0}^4a_k(x-x_i)^k.
\label{P4_function}
\end{equation}
The superscript `` $P_4$ '' stands for the 4th-degree polynomial. The unknown coefficients $a_k(k=0,1,2,3,4)$ are determined from the following constraint conditions,
\begin{equation}
\frac{1}{\Delta x}\int_{x_{j-\frac{1}{2}}}^{x_{j+\frac{1}{2}}}\mathcal{Q}_i^{P_4}(x)dx=\bar{q}_j\ \ (j=i,i\pm 1,i\pm 2),
\label{P4_condition}
\end{equation}
which require the volume average of the reconstruction function $\mathcal{Q}_i^{P_4}$ on cell $\Omega_j$ in the neighbors of the target cell $\Omega_i$ to be equal to the corresponding cell average of the solution $\bar{q}_j$. Consequently, the reconstructed boundary values can be obtained by
\begin{eqnarray}
\begin{aligned}
&q_{i+\frac{1}{2}}^{L,P_4}= \mathcal{Q}_i^{P_4}\left(x_{i+\frac{1}{2}}\right)=\frac{1}{60}(2\bar{q}_{i-2}-13\bar{q}_{i-1}+47\bar{q}_{i}+27\bar{q}_{i+1}-3\bar{q}_{i+2}), \\
&q_{i-\frac{1}{2}}^{R,P_4}= \mathcal{Q}_i^{P_4}\left(x_{i-\frac{1}{2}}\right)=\frac{1}{60}(2\bar{q}_{i+2}-13\bar{q}_{i+1}+47\bar{q}_{i}+27\bar{q}_{i-1}-3\bar{q}_{i-2}).
\label{P4_boundary}
\end{aligned}
\end{eqnarray}
For smooth solutions, the 4th-degree polynomial reconstruction results in a 5th-order convergence rate and provides a optimal approximation for smooth solution, which however tends to generate spurious oscillation due to the Gibbs phenomenon near discontinuous solutions. To suppress the oscillation, the 4th-degree polynomial is replaced by the THINC function  that is more suitable to approximate a discontinuous solution through the BVD selection algorithm.

\subsubsection{Candidate interpolant 2: THINC scheme}

THINC scheme was developed for capturing moving interfaces or free boundaries in multi-phase flows \cite{Xiao2005,Xiao2011,Ii2014,Xie2014,Xie2017}. Thanks to its  monotonicity and capability to mimic the step-like profile,  the THINC reconstruction is also well suited for computing compressible flows including discontinuous solutions with appealing stability and high resolution \cite{Deng2018,Deng2019,Deng2020,Tann2019,Tann2020,Cheng2021}. In the THINC scheme, the physical quantity $q$ is reconstructed by a hyperbolic tangent function as
\begin{eqnarray}
\mathcal{Q}_i^T(x)=
\bar{q}_{\rm min}+\frac{\Delta \bar{q}}{2}\left(1+\theta \tanh \left(\beta \left(X_i-d_i\right)\right)\right)
\label{THINC_function}
\end{eqnarray}
where
\begin{eqnarray*}
&&\bar{q}_{\rm min}={\rm min}(\bar{q}_{i-1}, \bar{q}_{i+1}), \ \Delta \bar{q}=|\bar{q}_{i+1}-\bar{q}_{i-1}|, \\
&&\theta={\rm sgn}(\bar{q}_{i+1}-\bar{q}_{i-1}), \ X_i=\frac{x-x_{i-\frac{1}{2}}}{x_{i+\frac{1}{2}}-x_{i-\frac{1}{2}}}.
\end{eqnarray*}
A superscript `` T '' stands for the THINC scheme. Since the THINC function  \eqref{THINC_function} is a monotonic function, it can be applied only if the monotonicity condition $(\bar{q}_i-\bar{q}_{i-1})(\bar{q}_{i+1}-\bar{q}_i)>0$ is satisfied. Or otherwise, the THINC function degrades to a piecewise constant reconstruction in practice, i.e. $\mathcal{Q}_i^T(x)=\bar{q}_i$.

Given the steepness parameter $\beta$, the only unknown parameter in \eqref{THINC_function} is the jump location (the center of the jump transition layer) $d_i$. We determine $d_i$ from the condition that the volume integrated average of the reconstruction function  $\mathcal{Q}_i^{T}$ over the target cell $\Omega_i$ is equal to the solution $\bar{q}_i$, which is written as
 \begin{equation}
\frac{1}{\Delta x}\int_{x_{i-\frac{1}{2}}}^{x_{i+\frac{1}{2}}}\mathcal{Q}_i^T(x)dx=\bar{q}_i.
\label{THINC_condition}
\end{equation}
Then, $d_i$ can be calculated by 
\begin{equation}
d_i=\frac{1}{2\beta}\ln{\frac{1-A}{1+A}},
\label{THINC_d}
\end{equation}
where $A=\frac{\exp (\alpha_i\beta)/\cosh(\beta)-1}{\tanh(\beta)}$, and the parameter $\alpha_i$, which implies the location of the discontinuity, is expressed as
\begin{eqnarray}
\alpha_i=\theta \left( 2\ \frac{\bar{q}_i-\bar{q}_{\rm min}+\epsilon}{\Delta \bar{q}+\epsilon}-1 \right).
\label{THINC_alpha}
\end{eqnarray}
The small number $\epsilon=10^{-20}$ is introduced to prevent arithmetic failure. Although we can obtain the boundary values of THINC scheme by calculating $d_i$ from \eqref{THINC_d} numerically and substituting it into function \eqref{THINC_function}, we directly compute the boundary values by the following formulae for computational efficiency,  
\begin{eqnarray}
\begin{aligned}
&q_{i+\frac{1}{2}}^{L,T}= \mathcal{Q}_i^{T}\left(x_{i+\frac{1}{2}}\right)=\bar{q}_{\rm min}+\frac{\Delta \bar{q}}{2}\left(1+\theta\frac{\tanh(\beta)+A}{1+A\tanh(\beta)}\right), \\
&q_{i-\frac{1}{2}}^{R,T}= \mathcal{Q}_i^{T}\left(x_{i-\frac{1}{2}}\right)=\bar{q}_{\rm min}+\frac{\Delta \bar{q}}{2}\left(1+\theta A\right).
\end{aligned}
\end{eqnarray}

The performance of THINC scheme depends on the value of the steepness  parameter $\beta$ as discussed in \cite{Deng2018,Tann2019}. When $\beta$ is set to 1.6, the reconstruction function becomes closer to a step-like profile and the discontinuous solution can be resolved within about four mesh cells \cite{Xiao2011}. On the other hand, the THINC scheme with $\beta=1.1$ leads to a smooth interpolation and has almost the same feature as the MUSCL scheme with the van Leer slope limiter \cite{VanLeer1979}, which has been demonstrated by the ADR (approximate dispersion relation) analysis \cite{Pirozzoli2006,Deng2018} and numerical results. We denote the THINC scheme with $\beta=1.1$ by THINC($\beta_s$) and that with $\beta=1.6$ by THINC($\beta_l$).  In the P$_4$T$_2$-BVD scheme, THINC($\beta_s$) is used to suppress numerical oscillation near the discontinuities, and  THINC($\beta_l$) is used to reduce the numerical dissipation in the vicinity of discontinuous solutions. 

\subsubsection{2-stage BVD algorithm}
As mentioned above, the selection of the reconstruction function is conducted via a  2-stage BVD algorithm. I.e., in the 1st stage, either the 4th-degree polynomial or THINC($\beta_s$) is selected (we denote the selected function as $\mathcal{Q}_i^{I}$), and in the 2nd stage, either the function selected in the 1st stage ($\mathcal{Q}_i^{I}$) or THINC($\beta_l$) is selected as the final reconstruction function.

We summarize the solution procedure of the 2-stage BVD algorithm for target cell $\Omega_i$ as follows.
\begin{enumerate}[(I)]

    \item \label{item1} The 1st stage
    \begin{enumerate}[(\ref{item1}-I)]
    
        \item Set the 4th-degree polynomial $\mathcal{Q}_i^{P_4}$ as the basic reconstruction function by 
        \begin{equation}
        \mathcal{Q}_i^{I}=\mathcal{Q}_i^{P_4}.
        \end{equation}
        
        \item Calculate the values of $TBV_i$ from the reconstruction of 4th-degree polynomial and THINC($\beta_s$) as
        \begin{eqnarray}
        \begin{aligned}
        &TBV_i^{P_4}=\left| q_{i-\frac{1}{2}}^{L,P_4}-q_{i-\frac{1}{2}}^{R,P_4}\right|+\left| q_{i+\frac{1}{2}}^{L,P_4}-q_{i+\frac{1}{2}}^{R,P_4}\right|, \\
        &TBV_i^{T_s}=\left| q_{i-\frac{1}{2}}^{L,T_s}-q_{i-\frac{1}{2}}^{R,T_s}\right|+\left| q_{i+\frac{1}{2}}^{L,T_s}-q_{i+\frac{1}{2}}^{R,T_s}\right|.
        \end{aligned}
        \end{eqnarray}
        
        \item Change the reconstruction functions for cells $\Omega_{i-1},\Omega_{i},\Omega_{i+1}$ to THINC($\beta_s$) according to the following BVD algorithm,
        \begin{eqnarray}
        \mathcal{Q}_j^{I}=\mathcal{Q}_j^{T_s}\ \ (j=i-1,i,i+1), \   {\rm if} \ \ TBV_i^{T_s}<TBV_i^{P_4}.
        \end{eqnarray}

    \end{enumerate}
    
    \item \label{item2} The 2nd stage
    \begin{enumerate}[(\ref{item2}-I)]
    
        \item Calculate the values of $TBV_i$ from the reconstruction of the function selected in the 1st stage and THINC($\beta_l$) as
        \begin{eqnarray}
        \begin{aligned}
        &TBV_i^{I}=\left| q_{i-\frac{1}{2}}^{L,I}-q_{i-\frac{1}{2}}^{R,I}\right|+\left| q_{i+\frac{1}{2}}^{L,I}-q_{i+\frac{1}{2}}^{R,I}\right|, \\
        &TBV_i^{T_l}=\left| q_{i-\frac{1}{2}}^{L,T_l}-q_{i-\frac{1}{2}}^{R,T_l}\right|+\left| q_{i+\frac{1}{2}}^{L,T_l}-q_{i+\frac{1}{2}}^{R,T_l}\right|.
        \end{aligned}
        \end{eqnarray}
        
        \item Determine the final reconstruction function for the target cell $\Omega_i$ by 
        \begin{eqnarray}
        \mathcal{Q}_i=
        \begin{cases}
        \mathcal{Q}_i^{T_l}, \ \ &{\rm if} \ \ TBV_i^{T_l}<TBV_i^{I}, \\
        \mathcal{Q}_i^{I}, &{\rm otherwise}.
        \end{cases}
        \end{eqnarray}
        
    \end{enumerate}
\end{enumerate}

In the 1st stage, the reconstruction functions in the neighboring cells of a discontinuity  are replaced by THINC($\beta_s$) to suppress numerical oscillations. In the 2nd stage, the reconstruction function in a cell containing the discontinuity is switched to THINC($\beta_l$) to reduce numerical dissipation. In both stages, the BVD principle that minimizes the reconstructed BV values is implemented, which not only effectively reduces numerical dissipation but also gives a reliable criterion to distinguish smooth and discontinuous solutions.

Shown in \cite{Deng2019}, the 2-stage BVD algorithm maintains the 5th-order convergence rate for smooth solutions, while suppresses spurious oscillations and gets stable numerical results for discontinuous solutions even without using the nonlinear limiting projection as those in the WENO-type high-resolution schemes. It is demonstrated that the P$_4$T$_2$-BVD scheme can retrieve the linear 4th-degree polynomial as the basic reconstruction function for smooth solutions, and is able to reproduce high-fidelity numerical results for both smooth and discontinuous solutions with superior solution quality in comparison with other existing high-order nonlinear shock-capturing schemes.

\section{Symmetry-preserving techniques}
\label{sec:Symmetry}

Even though numerical schemes are formulated to be mathematically exact, the causes of the breaking-symmetry might still stem from the mismatched rounding errors in the calculation of the physical quantities at symmetric positions \cite{Remacle2003,Fleischmann2019}. For example, the associativities of summation and multiplication for more than two components are generally lost in floating-point arithmetic as following,
\begin{equation}
(a+b)+c\neq a+(b+c),
\label{lack_of_associativity_sum}
\end{equation}
\begin{equation}
(a\times b)\times c\neq a\times(b\times c).
\label{lack_of_associativity_multiple}
\end{equation}
Thus, if the same formulae are computed with different orders of summation or multiplication, the numerical results might be different. Hence, the symmetricity in numerical solution is not automatically guaranteed, even the scheme is mathematically symmetric. In other words, the symmetry property cannot be completely preserved unless all processes of digit operations at symmetric positions are identical. In this section, we thoroughly examine and solve all possible symmetry-breaking elements in the numerical solution procedure of  the P$_4$T$_2$-BVD scheme, including  the transformation between conservative variables and characteristic variables, the spatial reconstruction, and the Riemann solver.

\subsection{The symmetry configurations in 2D simulations}

\begin{figure}[htbp]
    \centering
    \includegraphics[width=16.0cm]{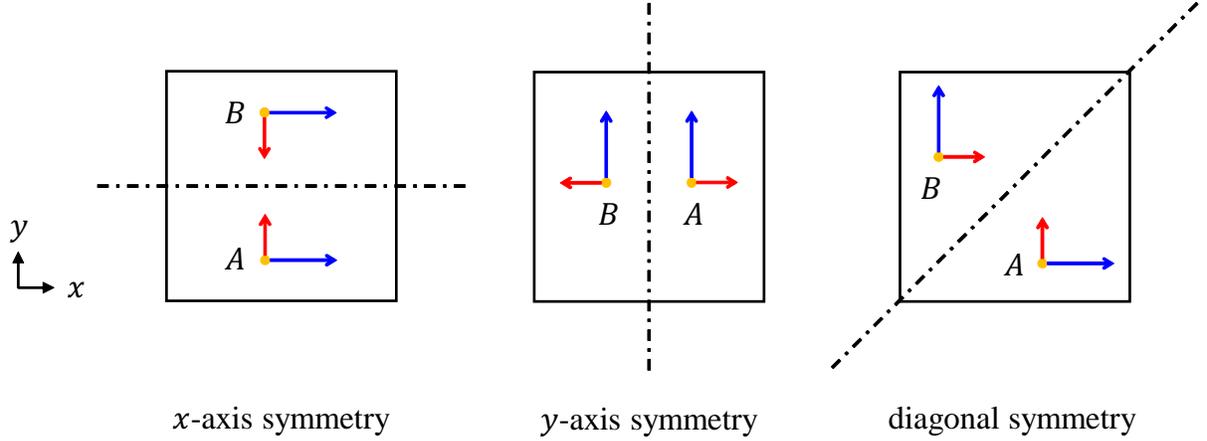}
    \caption{The three typical symmetry properties in 2D simulations.}
    \label{fig:symmetry_type}
\end{figure}

In this work, we consider three typical symmetrical configurations for Cartesian grid or structured grids in two dimensions, i.e.  $x$- or $y$-axis symmetry and diagonal symmetry, as illustrated in Fig. \ref{fig:symmetry_type}.
The square and the dashed line in Fig. \ref{fig:symmetry_type} show the computational domain and symmetry axis respectively. Marks $A$ and $B$ indicate the symmetric positions. Note that another type of diagonal symmetry property (the symmetry with respect to the  diagonal line cutting through Quadrants 2 and 4), which is not shown in Fig. \ref{fig:symmetry_type}, is automatically preserved if the three symmetry properties shown in Fig. \ref{fig:symmetry_type} are preserved. In 2D problems, preserving the three properties is the necessary and sufficient condition for eliminating the symmetry errors. 

The relationships of values of primitive variables $\rho,u,v,p$ at the positions $A$ and $B$ in each symmetry configuration are specified as follows,
\begin{align}
\label{prim_x-axis}
&\rho_A=\rho_B,\ u_A=u_B,\ v_A=-v_B,\ p_A=p_B\ \ \text{($x$-axis symmetry)}, \\
\label{prim_y-axis}
&\rho_A=\rho_B,\ u_A=-u_B,\ v_A=v_B,\ p_A=p_B\ \ \text{($y$-axis symmetry)}, \\
\label{prim_diag}
&\rho_A=\rho_B,\ u_A=v_B,\ v_A=u_B,\ p_A=p_B\ \ \text{(diagonal symmetry)},
\end{align}
which are referred to as the symmetry rules. 

The causes of the symmetry errors are investigated based on relationships \eqref{prim_x-axis}, \eqref{prim_y-axis} and \eqref{prim_diag}. Below, we use breve accent mark `` $\breve{\ }$ '' to highlight the symmetry-preserving formulations to distinguish them from the original ones.

\subsection{Local characteristic decomposition}

The reconstruction of the high-order schemes including the P$_4$T$_2$-BVD scheme is usually  performed in terms of the characteristic variables $\mathbf{\bar{W}}$ rather than  the  conservative variables $\mathbf{\bar{U}}$. So, the conservative variables need to be locally projected to the characteristic variables, and then the cell boundary values of reconstructed characteristic variables need to be projected back to conservative or primitive variables for use in Riemann solvers. The transformations between the conservative variables and characteristic variables are conducted by \eqref{trans-x} and  \eqref{trans-y}.

All the components in eigenvectors/matrices are calculated by density-based Roe averages of the reconstructed values at cell boundary, which reads, as an example in $x$-direction,
\begin{align}
q^{Roe}_{i+\frac{1}{2}}=\frac{\sqrt{{\rho}^L_{i+\frac{1}{2}}}{q}^L_{i+\frac{1}{2}}+\sqrt{{\rho}^R_{i+\frac{1}{2}}}{q}^R_{i+\frac{1}{2}}}{\sqrt{{\rho}^L_{i+\frac{1}{2}}}+\sqrt{{\rho}^R_{i+\frac{1}{2}}}},
\end{align}
where $q$ represents the physical quantity to be averaged, i.e., $u$, $v$, or $H$. 
Other formulations of the Roe average have been reported to affect the symmetry properties in \cite{Fleischmann2019}. 

The numerical formulations using the eigenvectors/matrices to transform conservative variables and characteristic variables might, unfortunately, cause symmetry errors unless special manipulations are applied. 

First, let us consider the symmetry with respect to $y$-axis. In this case, the symmetry property is relevant only to the numerical formulations in $x$-direction, so we examine  the projection in $x$-direction. 

The characteristic variables $\mathbf{\bar{W}}_x$ at symmetric positions $A$ and $B$ in the middle panel of   Fig. \ref{fig:symmetry_type} are obtained by $\mathbf{\bar{W}}_{xA}=\mathbf{L}_{xA}\cdot\mathbf{\bar{U}}$ and $\mathbf{\bar{W}}_{xB}=\mathbf{L}_{xB}\cdot\mathbf{\bar{U}}$ respectively. We write explicitly  $\mathbf{\bar{W}}_{xA}$ and $\mathbf{\bar{W}}_{xB}$ using the symmetry relationship with respect to $y$-axis Eq. \eqref{prim_y-axis} as follows,
\begin{align}
\begin{pmatrix}
\bar{w}^{(u-c)} \\ \bar{w}^{(u)} \\ \bar{w}^{(u+c)} \\ \bar{w}^{(u\perp)}
\end{pmatrix}
_A&=
\begin{pmatrix}
\frac{1}{2}\left(b_1+\frac{u}{c}\right) \\ 1-b_1 \\ \frac{1}{2}\left(b_1-\frac{u}{c}\right) \\ -v
\end{pmatrix}
_A\bar{\rho}_A+
\begin{pmatrix}
-\frac{1}{2}\left(\frac{1}{c}+b_2u\right) \\ b_2u \\ \frac{1}{2}\left(\frac{1}{c}-b_2u\right) \\ 0
\end{pmatrix}
_A\bar{\rho u}_A+
\begin{pmatrix}
-\frac{1}{2}b_2v \\ b_2v \\ -\frac{1}{2}b_2v \\ 1
\end{pmatrix}
_A\bar{\rho v}_A+
\begin{pmatrix}
\frac{1}{2}b_2 \\ -b_2 \\ \frac{1}{2}b_2 \\ 0
\end{pmatrix}
_A\bar{E}_A, \nonumber \\
\begin{pmatrix}
\bar{w}^{(u-c)} \\ \bar{w}^{(u)} \\ \bar{w}^{(u+c)} \\ \bar{w}^{(u\perp)}
\end{pmatrix}
_B&=
\begin{pmatrix}
\frac{1}{2}\left(b_1+\frac{u}{c}\right) \\ 1-b_1 \\ \frac{1}{2}\left(b_1-\frac{u}{c}\right) \\ -v
\end{pmatrix}
_B\bar{\rho}_B+
\begin{pmatrix}
-\frac{1}{2}\left(\frac{1}{c}+b_2u\right) \\ b_2u \\ \frac{1}{2}\left(\frac{1}{c}-b_2u\right) \\ 0
\end{pmatrix}
_B\bar{\rho u}_B+
\begin{pmatrix}
-\frac{1}{2}b_2v \\ b_2v \\ -\frac{1}{2}b_2v \\ 1
\end{pmatrix}
_B\bar{\rho v}_B+
\begin{pmatrix}
\frac{1}{2}b_2 \\ -b_2 \\ \frac{1}{2}b_2 \\ 0
\end{pmatrix}
_B\bar{E}_B \nonumber \\
&=
\begin{pmatrix}
\frac{1}{2}\left(b_1-\frac{u}{c}\right) \\ 1-b_1 \\ \frac{1}{2}\left(b_1+\frac{u}{c}\right) \\ -v
\end{pmatrix}
_A\bar{\rho}_A+
\begin{pmatrix}
-\frac{1}{2}\left(\frac{1}{c}-b_2u\right) \\ -b_2u \\ \frac{1}{2}\left(\frac{1}{c}+b_2u\right) \\ 0
\end{pmatrix}
_A(-\bar{\rho u})_A+
\begin{pmatrix}
-\frac{1}{2}b_2v \\ b_2v \\ -\frac{1}{2}b_2v \\ 1
\end{pmatrix}
_A\bar{\rho v}_A+
\begin{pmatrix}
\frac{1}{2}b_2 \\ -b_2 \\ \frac{1}{2}b_2 \\ 0
\end{pmatrix}
_A\bar{E}_A \nonumber \\
&=
\begin{pmatrix}
\bar{w}^{(u+c)} \\ \bar{w}^{(u)} \\ \bar{w}^{(u-c)} \\ \bar{w}^{(u\perp)}
\end{pmatrix}
_A.
\end{align}
\begin{figure}[htbp]
\centering
\includegraphics[width=12.0cm]{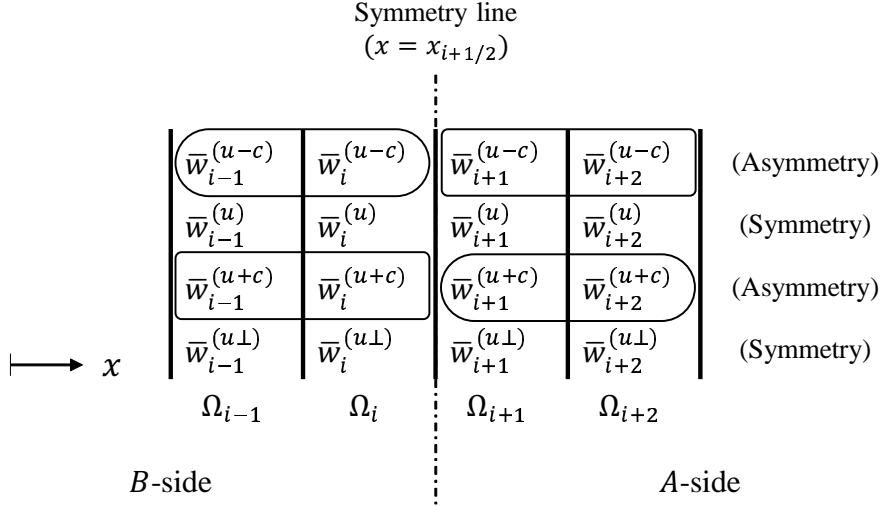}
\caption{The distributions of the characteristic variables  $\mathbf{\bar{W}}_x$. The variables boxed by squares and the variables boxed by ellipses on the two sides of the symmetry line $x=x_{i+1/2}$ interchange their positions.} 
\label{fig:dist_Wx}
\end{figure}

It is observed from this expressions that the values of $\bar{w}^{(u-c)}$ and $\bar{w}^{(u+c)}$ are interchanged while the values of $\bar{w}^{(u)}$ and $\bar{w}^{(u\perp)}$  maintain the same at the symmetric positions $A$ and $B$. In other words, the symmetry of $\bar{w}^{(u-c)}$ and $\bar{w}^{(u+c)}$ are broken although the conservative variables $\mathbf{\bar{U}}$ are exactly symmetrical. The distributions of the characteristic variables  $\mathbf{\bar{W}}_x$  are illustrated in Fig. \ref{fig:dist_Wx}. It is seen that the cell boundary values of $w^{(u-c)}$ and $w^{(u+c)}$ obtained by reconstruction are not the same about symmetric line  $x=x_{i+1/2}$. Nevertheless, the symmetry in conservative variables can be recovered if we retrieve the conservative variables at position $A$ and $B$  by  $\mathbf{U}_A=\mathbf{R}_{xA}\cdot\mathbf{W}_{xA}$ and $\mathbf{U}_B=\mathbf{R}_{xB}\cdot\mathbf{W}_{xB}$ respectively.

Recall that 
 \begin{align}
\begin{pmatrix}
\rho \\ \rho u \\ \rho v \\ E
\end{pmatrix}
_A&=
\begin{pmatrix}
1 \\ u-c \\ v \\ H-uc
\end{pmatrix}
_Aw^{(u-c)}_A+
\begin{pmatrix}
1 \\ u \\ v \\ \frac{u^2+v^2}{2}
\end{pmatrix}
_Aw^{(u)}_A+
\begin{pmatrix}
1 \\ u+c \\ v \\ H+uc
\end{pmatrix}
_Aw^{(u+c)}_A+
\begin{pmatrix}
0 \\ 0 \\ 1 \\ v
\end{pmatrix}
_Aw^{(u\perp)}_A, \nonumber \\
\end{align}It is straightforward to prove the following relationship with symmetry rule  \eqref{prim_y-axis}, 
\begin{align}
\begin{pmatrix}
\rho \\ \rho u \\ \rho v \\ E
\end{pmatrix}
_B&=
\begin{pmatrix}
1 \\ u-c \\ v \\ H-uc
\end{pmatrix}
_Bw^{(u-c)}_B+
\begin{pmatrix}
1 \\ u \\ v \\ \frac{u^2+v^2}{2}
\end{pmatrix}
_Bw^{(u)}_B+
\begin{pmatrix}
1 \\ u+c \\ v \\ H+uc
\end{pmatrix}
_Bw^{(u+c)}_B+
\begin{pmatrix}
0 \\ 0 \\ 1 \\ v
\end{pmatrix}
_Bw^{(u\perp)}_B \nonumber \\ &=
\begin{pmatrix}
1 \\ {-(u+c)} \\ v \\ H+uc
\end{pmatrix}
_Aw^{(u+c)}_A+
\begin{pmatrix}
1 \\ -u \\ v \\ \frac{u^2+v^2}{2}
\end{pmatrix}
_Aw^{(u)}_A+
\begin{pmatrix}
1 \\ {-(u-c)} \\ v \\ H-uc
\end{pmatrix}
_Aw^{(u-c)}_A+
\begin{pmatrix}
0 \\ 0 \\ 1 \\ v
\end{pmatrix}
_Aw^{(u\perp)}_A \nonumber \\ 
&{=}\begin{pmatrix}
\rho \\ -\rho u \\ \rho v \\ E
\end{pmatrix}
_A, 
\label{WB_to_UB_y-axis}
\end{align}
which states that the projection between conservative variables and characteristic variables does not break the symmetry in theory if there is no any numerical error involved. 

However, Eq. \eqref{WB_to_UB_y-axis} does not necessarily hold in numerical processing due to the floating-point arithmetic operation, which might consequently lead to asymmetry in numerical solution. As shown before in \eqref{lack_of_associativity_sum} and \eqref{lack_of_associativity_multiple}, the floating-point arithmetic does not ensure the  rigorous associativity of algebraic operations because of the round-off error. As an example, we examine the calculation procedures of $\rho_A$ and $\rho_B$, 
\begin{align}
\rho_A&=(\underset{\text{a}}{\uwave{w^{(u-c)}_A}}+\underset{\text{b}}{\uwave{w^{(u)}_A}})+\underset{\text{c}}{\uwave{w^{(u+c)}_A}}, \nonumber \\
\rho_B&=(\underset{\text{c}}{\uwave{w^{(u-c)}_B}}+\underset{\text{b}}{\uwave{w^{(u)}_B}})+\underset{\text{a}}{\uwave{w^{(u+c)}_B}} \nonumber \\
&=(\underset{\text{c}}{\uwave{w^{(u+c)}_A}}+\underset{\text{b}}{\uwave{w^{(u)}_A}})+\underset{\text{a}}{\uwave{w^{(u-c)}_A}} \nonumber \\
&\neq\rho_A.
\end{align}
Here, the terms labeled with the same letter have the same values.
It observed that they might arrive at different values because of using different orders of summations. 
This observation applies also to other conservative variables, $\rho u,\rho v$ and $E$,  and eventually results in the symmetry errors. 

A remedy to this problem can be regrouping and rearranging the order of the operations in the transformation calculation. I.e., sum up first the variables corresponding to the eigenvalues $(u-c)$ and $(u+c)$ in $x$-direction and $(v-c)$ and $(v+c)$ in $y$-direction first as follows, 
\begin{eqnarray}
\begin{aligned}
&\breve{\mathbf{U}}=\pmb{\bigl(}\vec{r}^{(u-c)}w^{(u-c)}+\vec{r}^{(u+c)}w^{(u+c)}\pmb{\bigr)}+\vec{r}^{(u)}w^{(u)}+\vec{r}^{(u\perp)}w^{(u\perp)}\ \ \text{($x$-direction)}, \\
&\breve{\mathbf{U}}=\pmb{\bigl(}\vec{r}^{(v-c)}w^{(v-c)}+\vec{r}^{(v+c)}w^{(v+c)}\pmb{\bigr)}+\vec{r}^{(v)}w^{(v)}+\vec{r}^{(v\perp)}w^{(v\perp)}\ \ \text{($y$-direction)}.
\end{aligned}
\label{solution_char_y-axis}
\end{eqnarray}
Implementing this method requires to modify the calculation procedure for matrix-vector product. 
Alternatively,  the same can be achieved by rearranging the order of the eigenvalues and eigenvectors in matrices $\mathbf{L}$ and $\mathbf{R}$. To this end,  we take  $(u-c)$ and $(u+c)$ as the first two eigenvalues, and rearrange the eigenvectors to formulate the so-called ``symmetry-preserving eigenvectors matrices'',  $\breve{\mathbf{L}}$ and $\breve{\mathbf{R}}$ as follows, 
\begin{align}
\begin{aligned}
&\breve{\mathbf{L}}_x=
\begin{pmatrix}
\boxed{\vec{l}^{(u-c)}} \\ \boxed{\vec{l}^{(u+c)}} \\ \vec{l}^{(u)} \\ \vec{l}^{(u\perp)}
\end{pmatrix}
,\ \breve{\mathbf{L}}_y=
\begin{pmatrix}
\boxed{\vec{l}^{(v-c)}} \\ \boxed{\vec{l}^{(v+c)}} \\ \vec{l}^{(v)} \\ \vec{l}^{(v\perp)}
\end{pmatrix}
,\\&\breve{\mathbf{R}}_x=
\begin{pmatrix}
\boxed{\vec{r}^{(u-c)}} & \boxed{\vec{r}^{(u+c)}} & \vec{r}^{(u)} & \vec{r}^{(u\perp)}
\end{pmatrix}
,\\ &\breve{\mathbf{R}}_y=
\begin{pmatrix}
\boxed{\vec{r}^{(v-c)}} & \boxed{\vec{r}^{(v+c)}} & \vec{r}^{(v)} & \vec{r}^{(v\perp)}
\end{pmatrix},
\label{eigen_sym}
\end{aligned}
\end{align}
where the boxed eigen vectors are placed at the first two positions to make sure that the summation of them is computed first.  

Hence, the symmetry property regarding $x$- and $y$-axis can be  preserved by either adjusting the order of the summation of matrix-vector product as shown in \eqref{solution_char_y-axis} or re-arranging the eigen matrices as \eqref{eigen_sym}. 


We now analyse the diagonal symmetry shown in the right panel of Fig.\ref{fig:symmetry_type}. Consider two symmetric points $A$ and $B$, we write  the characteristic variables $\mathbf{\bar{W}}_{xA}$ at point $A$ as $\mathbf{\bar{W}}_{xA}=\mathbf{L}_{xA}\cdot\mathbf{\bar{U}}$
\begin{align}
\begin{pmatrix}
\bar{w}^{(u-c)} \\ \bar{w}^{(u)} \\ \bar{w}^{(u+c)} \\ \bar{w}^{(u\perp)}
\end{pmatrix}
_A&=
\begin{pmatrix}
\frac{1}{2}\left(b_1+\frac{u}{c}\right) \\ 1-b_1 \\ \frac{1}{2}\left(b_1-\frac{u}{c}\right) \\ -v
\end{pmatrix}
_A\bar{\rho}_A+
\begin{pmatrix}
-\frac{1}{2}\left(\frac{1}{c}+b_2u\right) \\ b_2u \\ \frac{1}{2}\left(\frac{1}{c}-b_2u\right) \\ 0
\end{pmatrix}
_A\bar{\rho u}_A+
\begin{pmatrix}
-\frac{1}{2}b_2v \\ b_2v \\ -\frac{1}{2}b_2v \\ 1
\end{pmatrix}
_A\bar{\rho v}_A+
\begin{pmatrix}
\frac{1}{2}b_2 \\ -b_2 \\ \frac{1}{2}b_2 \\ 0
\end{pmatrix}
_A\bar{E}_A. 
\end{align}
Using the symmetry rule \eqref{prim_diag}, we can get the characteristic variables  $\mathbf{\bar{W}}_{yB}=\mathbf{L}_{yB}\cdot\mathbf{\bar{U}}$ at point $B$ as, 
\begin{align}
\begin{pmatrix}
\bar{w}^{(v-c)} \\ \bar{w}^{(v)} \\ \bar{w}^{(v+c)} \\ \bar{w}^{(v\perp)}
\end{pmatrix}
_B&=
\begin{pmatrix}
\frac{1}{2}\left(b_1+\frac{v}{c}\right) \\ 1-b_1 \\ \frac{1}{2}\left(b_1-\frac{v}{c}\right) \\ -u
\end{pmatrix}
_B\bar{\rho}_B+
\begin{pmatrix}
-\frac{1}{2}b_2u \\ b_2u \\ -\frac{1}{2}b_2u \\ 1
\end{pmatrix}
_B\bar{\rho u}_B+
\begin{pmatrix}
-\frac{1}{2}\left(\frac{1}{c}+b_2v\right) \\ b_2v \\ \frac{1}{2}\left(\frac{1}{c}-b_2v\right) \\ 0
\end{pmatrix}
_B\bar{\rho v}_B+
\begin{pmatrix}
\frac{1}{2}b_2 \\ -b_2 \\ \frac{1}{2}b_2 \\ 0
\end{pmatrix}
_B\bar{E}_B \nonumber \\
&=
\begin{pmatrix}
\frac{1}{2}\left(b_1+\frac{u}{c}\right) \\ 1-b_1 \\ \frac{1}{2}\left(b_1-\frac{u}{c}\right) \\ -v
\end{pmatrix}
_A\bar{\rho}_A+
\begin{pmatrix}
-\frac{1}{2}b_2v \\ b_2v \\ -\frac{1}{2}b_2v \\ 1
\end{pmatrix}
_A\bar{\rho v}_A+
\begin{pmatrix}
-\frac{1}{2}\left(\frac{1}{c}+b_2u\right) \\ b_2u \\ \frac{1}{2}\left(\frac{1}{c}-b_2u\right) \\ 0
\end{pmatrix}
_A\bar{\rho u}_A+
\begin{pmatrix}
\frac{1}{2}b_2 \\ -b_2 \\ \frac{1}{2}b_2 \\ 0
\end{pmatrix}
_A\bar{E}_A \nonumber \\ &{=}
\begin{pmatrix}
\bar{w}^{(u-c)} \\ \bar{w}^{(u)} \\ \bar{w}^{(u+c)} \\ \bar{w}^{(u\perp)}
\end{pmatrix}
_A.
\label{UB_to_WB_diag}
\end{align}
It indicates that the symmetry property of the characteristic variables, $\mathbf{\bar{W}}_{xA}=\mathbf{\bar{W}}_{yB}$, is mathematically true. However, this symmetry is not guaranteed due to the loss of summation associativity  in floating-point arithmetic. We follow the  technique introduced in  \cite{Fleischmann2019}  and add the brackets to  adjust the sequence of the summation, 
\begin{eqnarray}
\begin{aligned}
\mathbf{\breve{\bar{W}}}_x=\vec{l}_{x1}\bar{\rho}+\pmb{\Bigl(}\vec{l}_{x2}\bar{\rho u}+\vec{l}_{x3}\bar{\rho v}\pmb{\Bigr)}+\vec{l}_{x4}\bar{E},
\end{aligned}
\label{solution_char_diag_x}\\
\begin{aligned}
\mathbf{\breve{\bar{W}}}_y=\vec{l}_{y1}\bar{\rho}+\pmb{\Bigl(}\vec{l}_{y2}\bar{\rho u}+\vec{l}_{y3}\bar{\rho v}\pmb{\Bigr)}+\vec{l}_{y4}\bar{E},
\end{aligned}
\label{solution_char_diag_y}
\end{eqnarray}
where $\vec{l}_{xm}$ and  $\vec{l}_{ym}$,  $m=1,2,3,4$, are the column vectors of  $\mathbf{L}_{x}$ and  $\mathbf{L}_{y}$ respectively in \eqref{eigen_L}.  It is noted that the transformation from the reconstructed characteristic variables to the conservative variables, i.e.  $\mathbf{U}=\mathbf{R}_{xA}\cdot\mathbf{W}_{xA}$ and $\mathbf{U}=\mathbf{R}_{yB}\cdot\mathbf{W}_{yB}$ in $x$- and $y$-directions respectively, don't break the symmetry property. We prove it in the appendix.

To conclude this subsection regarding the transformation between conservative variables and characteristic variables, we state the two crucial points to enforce the solution symmetry, 1) to add a bracket and modify the summation order as \eqref{solution_char_y-axis} or arrange the eigenvectors matrices as \eqref{eigen_sym} to ensure the axis-symmetry property, and 2) to add a bracket to \eqref{solution_char_diag_x} and \eqref{solution_char_diag_y} for diagonal symmetry-preserving.

\subsection{Reconstruction}

In order to avoid the symmetry errors in the reconstruction step, it is necessary to perform the reconstruction symmetrically at the level of the floating-point arithmetic. Specifically, the reconstructed values at cell boundaries have to exactly follow the symmetric rules  \eqref{prim_x-axis}-\eqref{prim_diag} at the symmetric positions. 

\begin{table}[htbp]
    \centering
    \begin{tabular}{cccc}\hline
        & \shortstack{$x$-axis\\ symmetry} & \shortstack{$y$-axis\\ symmetry} & \shortstack{diagonal\\ symmetry} \\ \hline
        \multirow{7}{*}{\shortstack{conservative\\ and\\ primitive\\ variables}}
        & $\rho_A=\rho_B$ & $\rho_A=\rho_B$ & $\rho_A=\rho_B$ \\
        & $\rho u_A=\rho u_B$ & $\rho u_A=-\rho u_B$ & $\rho u_A=\rho v_B$ \\
        & $\rho v_A=-\rho v_B$ & $\rho v_A=\rho v_B$ & $\rho v_A=\rho u_B$ \\
        & $E_A=E_B$ & $E_A=E_B$ & $E_A=E_B$ \\
        & $u_A=u_B$ & $u_A=-u_B$ & $u_A=v_B$ \\
        & $v_A=-v_B$ & $v_A=v_B$ & $v_A=u_B$ \\
        & $p_A=p_B$ & $p_A=p_B$ & $p_A=p_B$ \\ \hline
        \multirow{8}{*}{\shortstack{characteristic\\ variables}} 
        & $w^{(u-c)}_A=w^{(u-c)}_B$ & $w^{(u-c)}_A=w^{(u+c)}_B$ & $w^{(u-c)}_A=w^{(v-c)}_B$ \\
        & $w^{(u)}_A=w^{(u)}_B$ & $w^{(u)}_A=w^{(u)}_B$ & $w^{(u)}_A=w^{(v)}_B$ \\
        & $w^{(u+c)}_A=w^{(u+c)}_B$ & $w^{(u+c)}_A=w^{(u-c)}_B$ & $w^{(u+c)}_A=w^{(v+c)}_B$ \\
        & $w^{(u\perp)}_A=w^{(u\perp)}_B$ & $w^{(u\perp)}_A=w^{(u\perp)}_B$ & $w^{(u\perp)}_A=w^{(v\perp)}_B$ \\
        & $w^{(v-c)}_A=w^{(v+c)}_B$ & $w^{(v-c)}_A=w^{(v-c)}_B$ & $w^{(v-c)}_A=w^{(u-c)}_B$ \\
        & $w^{(v)}_A=w^{(v)}_B$ & $w^{(v)}_A=w^{(v)}_B$ & $w^{(v)}_A=w^{(u)}_B$ \\
        & $w^{(v+c)}_A=w^{(v-c)}_B$ & $w^{(v+c)}_A=w^{(v+c)}_B$ & $w^{(v+c)}_A=w^{(u+c)}_B$ \\
        & $w^{(v\perp)}_A=w^{(v\perp)}_B$ & $w^{(v\perp)}_A=w^{(v\perp)}_B$ & $w^{(v\perp)}_A=w^{(u\perp)}_B$ \\ \hline
    \end{tabular}
    \caption{Summary of the symmetry relationships of the conservative, primitive and characteristic variables.}
    \label{tab:relationships_summary}
\end{table}

We summarize the required relations of reconstructed variables for different types of symmetry in table \ref{tab:relationships_summary}. 
It is found that the stencils and the sequences of calculations for the physical quantities are the same at the diagonally symmetric positions. Thus, the reconstructions do not break the diagonal symmetry because the  formulations used for reconstructions in $x$- and $y$-directions are identical as implemented in this work. However, the reconstructed cell boundary values at ($x$- or $y$-) axis-symmetric positions do not automatically guarantee the symmetric rules, and special attention must be paid to ensure the exact symmetry in numerical solutions. In present work, we propose conditions to conduct symmetric reconstructions that are derived by leveraging two special reconstruction techniques, stencil-flipping (SF) reconstruction and  sign-inversion (SI) reconstruction. 

\subsubsection{Symmetry enforcement using SF and SI reconstructions}

Without losing generality, consider a spatial reconstruction for variable $q(x)$ in $x$-direction, we denote the piecewise reconstruction function $\mathcal{Q}_i(x)$ for cell $\Omega_{i}$ by $ \mathcal{Q}_i(\{\bar{q}_{j} \}, x)$, where $\{\bar{q}_{j} \}$ stands for the union of the cell-average values in the neighboring cells $\{ \Omega_{j} \}, \ j=i-i_l,\cdots, i-1,i,i+1,\cdots,i+i_r$. 

The SF reconstruction function for  $\Omega_{i}$  is defined by flipping the stencil as 
\begin{equation}
\mathcal{Q}^{SF}_i(x)\equiv \mathcal{Q}_i(\{\bar{q}_{j'} \}, x),  
\label{sf-recon}
\end{equation}
where the order of cells $\{ \Omega_{j'} \}$ is reversed as, $\ j'=i-i_r,\cdots, i-1,i,i+1,\cdots,i+i_l$. 

The SI reconstruction function for  $\Omega_{i}$  is defined as a reconstruction for the cell-average values with opposite sign, 
\begin{equation}
\mathcal{Q}^{SI}_i(x)\equiv \mathcal{Q}_i(\{-\bar{q}_{j} \}, x). 
\end{equation}

Given the SF function defined in \eqref{sf-recon}, the symmetric reconstruction in $x$-direction, regarding the $y$-axis symmetry, requires the following relations to be satisfied,  
\begin{equation}
\label{cond_flip}
\mathcal{Q}_i(x_{i+\frac{1}{2}})=\mathcal{Q}^{SF}_i(x_{i-\frac{1}{2}}) \ \ \text{and} \ \ \mathcal{Q}_i(x_{i-\frac{1}{2}})=\mathcal{Q}^{SF}_i(x_{i+\frac{1}{2}}),  
\end{equation}
which applies to the reconstructions of all variables. 

The SI reconstruction needs to fulfill the following relation as an additional requisite to ensure the exact symmetry with respect to axis for the THINC function discussed later,  

\begin{equation}
\label{cond_invers}
\mathcal{Q}_i(x_{i+\frac{1}{2}})=-\mathcal{Q}^{SI}_i(x_{i+\frac{1}{2}}) \ \ \text{and} \ \ \mathcal{Q}_i(x_{i-\frac{1}{2}})=-\mathcal{Q}^{SI}_i(x_{i-\frac{1}{2}}).   
\end{equation}
We also note that for reconstructions of the conservative or primitive variables $\rho u$ or $u$,  which has opposite signs at the symmetric positions about $y$-axis, condition \eqref{cond_invers} is required for symmetry. It applies analogously to the reconstruction in $y$-direction for the symmetry about $x$-axis as well. 

It is found that the numerical formulations for reconstructions are not automatically in line with conditions \eqref{cond_flip} and \eqref{cond_invers}, and hence might cause the symmetry-breaking in numerical solutions. Next, we examine specifically the reconstruction functions used in the  P$_4$T$_2$-BVD method.

\subsubsection{Symmetry property of 4th-degree polynomial function}
The cell-boundary values of the 4th-degree polynomial reconstruction are computed as  linear combinations of $\bar{q}_j\ (j=i,i\pm1,i\pm2)$, i.e. the summation of five values. As discussed above, due to the invalid associativity of summation in the floating-point arithmetic, the orders of the summations might cause symmetry-breaking. Summing up the values in a natural order with the cell index gradually increased, leads to the following formulae to calculate cell-boundary values, which unfortunately does not preserve the symmetry, 
\begin{align}
\begin{aligned}
&q_{i+\frac{1}{2}}^{L,P_4}=\frac{1}{60}(\underset{\text{a}}{\uwave{2\bar{q}_{i-2}}}-\underset{\text{b}}{\uwave{13\bar{q}_{i-1}}}+\underset{\text{c}}{\uwave{47\bar{q}_{i}}}+\underset{\text{d}}{\uwave{27\bar{q}_{i+1}}}-\underset{\text{e}}{\uwave{3\bar{q}_{i+2}}}), \\
&q_{i-\frac{1}{2}}^{R,P_4}=\frac{1}{60}(-\underset{\text{e}}{\uwave{3\bar{q}_{i-2}}}+\underset{\text{d}}{\uwave{27\bar{q}_{i-1}}}+\underset{\text{c}}{\uwave{47\bar{q}_{i}}}-\underset{\text{b}}{\uwave{13\bar{q}_{i+1}}}+\underset{\text{a}}{\uwave{2\bar{q}_{i+2}}}),
\label{P4_boundary_asym}
\end{aligned}
\end{align}
where we label the five terms from ``a'' to ``e''. It is clear that \eqref{P4_boundary_asym} does not meet the symmetric condition   \eqref{cond_flip}. Applying \eqref{cond_flip} leads to another arrangement of the terms that matches position of the terms of the same label, which equivalently to flipping the stencil, i.e. $\bar{q}_{i-2}$ and $\bar{q}_{i+2}$, $\bar{q}_{i-1}$ and $\bar{q}_{i+1}$ are interchanged respectively. Hence, the order of the summation of $q_{i-\frac{1}{2}}^{R,P_4}$ are adjusted to that of $q_{i+\frac{1}{2}}^{L,P_4}$ as follows,
\begin{eqnarray}
\begin{aligned}
&\breve{q}_{i+\frac{1}{2}}^{L,P_4}=\frac{1}{60}(\underset{\text{a}}{\uwave{2\bar{q}_{i-2}}}-\underset{\text{b}}{\uwave{13\bar{q}_{i-1}}}+\underset{\text{c}}{\uwave{47\bar{q}_{i}}}+\underset{\text{d}}{\uwave{27\bar{q}_{i+1}}}-\underset{\text{e}}{\uwave{3\bar{q}_{i+2}}}), \\
&\breve{q}_{i-\frac{1}{2}}^{R,P_4}=\frac{1}{60}(\underset{\text{a}}{\uwave{2\bar{q}_{i+2}}}-\underset{\text{b}}{\uwave{13\bar{q}_{i+1}}}+\underset{\text{c}}{\uwave{47\bar{q}_{i}}}+\underset{\text{d}}{\uwave{27\bar{q}_{i-1}}}-\underset{\text{e}}{\uwave{3\bar{q}_{i-2}}}),
\label{P4_boundary_sym}
\end{aligned}
\end{eqnarray}
which exactly preserves the symmetry in numerical solutions.  In sense of eliminating symmetry errors in summation operation, the following formulae are effective as well by explicitly adding  brackets to group the values of the cells symmetric about $\Omega_{i}$.  
\begin{eqnarray}
\begin{aligned}
&\breve{q}_{i+\frac{1}{2}}^{L,P_4}=\frac{1}{60}((\underset{\text{a}}{\uwave{2\bar{q}_{i-2}}}-\underset{\text{e}}{\uwave{3\bar{q}_{i+2}}})+(-\underset{\text{b}}{\uwave{13\bar{q}_{i-1}}}+\underset{\text{d}}{\uwave{27\bar{q}_{i+1}}})+\underset{\text{c}}{\uwave{47\bar{q}_{i}}}), \\
&\breve{q}_{i-\frac{1}{2}}^{R,P_4}=\frac{1}{60}((-\underset{\text{e}}{\uwave{3\bar{q}_{i-2}}}+\underset{\text{a}}{\uwave{2\bar{q}_{i+2}}})+(\underset{\text{d}}{\uwave{27\bar{q}_{i-1}}}-\underset{\text{b}}{\uwave{13\bar{q}_{i+1}}})+\underset{\text{c}}{\uwave{47\bar{q}_{i}}}).
\label{P4_boundary_sym_lr}
\end{aligned}
\end{eqnarray}

It is noted that the SI reconstruction does not make anything different regarding symmetry property in the 4th-degree polynomial reconstruction because \eqref{cond_invers} always holds no matter if $\bar{q}_j\ (j=i,i\pm1,i\pm2)$ are changed to opposite sign or not. This conclusion applies to all linear schemes that use polynomials with constant coefficients for reconstruction. 

\subsubsection{A new formulation of THINC function for symmetry-preserving}

Here, we propose a new formulation of the THINC function to completely remove the symmetry errors. This formulation is devised
 to fulfill conditions \eqref{cond_flip} and \eqref{cond_invers}. Instead of \eqref{THINC_function}, the new THINC function is written as 

\begin{figure}[htbp]
    \centering
    \includegraphics[width=16.0cm]{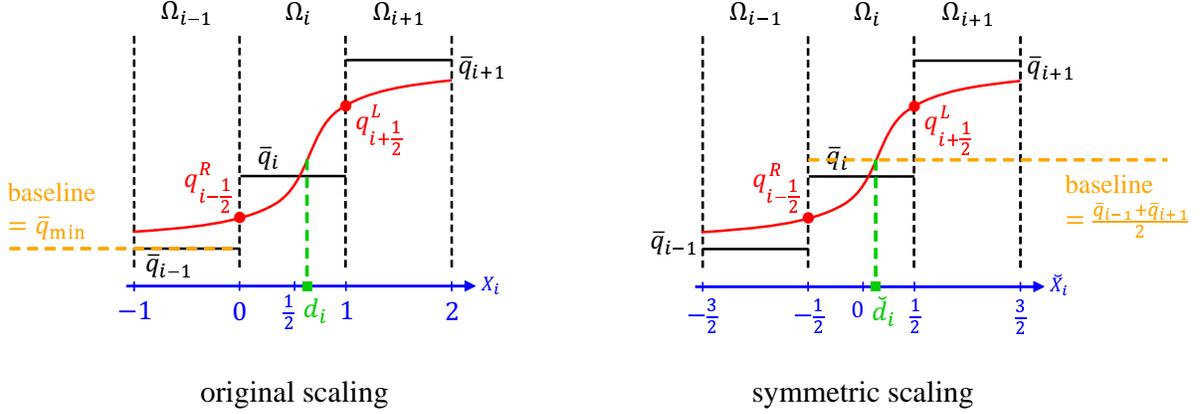}
    \caption{The descriptions of the scaling of original THINC and symmetric THINC function.}
    \label{fig:THINC_scaling}
\end{figure}

\begin{equation}
\breve{\mathcal{Q}}_i^{T}(x)={q}_{a}+{q}_{d}\tanh \left(\beta \left(\breve{X}_i-\breve{d}_i\right)\right),
\label{symTHINC_function}
\end{equation}
where
\begin{eqnarray*}
&&{q}_{a}=\frac{\bar{q}_{i+1}+\bar{q}_{i-1}}{2},\ {q}_{d}=\frac{\bar{q}_{i+1}-\bar{q}_{i-1}}{2}, \\
&&\breve{X}_i=X_i-\frac{1}{2}=\frac{x-\left( x_{i+\frac{1}{2}}+x_{i-\frac{1}{2}}\right) /2}{x_{i+\frac{1}{2}}-x_{i-\frac{1}{2}}}.
\end{eqnarray*}
Same as the original THINC scheme, function \eqref{symTHINC_function} is applied only if the monotonicity condition $(\bar{q}_i-\bar{q}_{i-1})(\bar{q}_{i+1}-\bar{q}_i)>10^{-20}$ is satisfied. Otherwise, the THINC function degrades to a piecewise constant function, i.e. $\breve{\mathcal{Q}}_i^T(x)=\bar{q}_i$. The jump location (center of the transition layer) $\breve{d}_i$ is obtained by the following condition, same as the original THINC scheme,
\begin{equation}
\frac{1}{\Delta x}\int_{x_{i-\frac{1}{2}}}^{x_{i+\frac{1}{2}}}\breve{\mathcal{Q}}_i^T(x)dx=\bar{q}_i.
\label{symTHINC_condition}
\end{equation}
Thus, $\breve{d}_i$ can be derived as
\begin{equation}
\breve{d}_i=\frac{1}{2\beta}\ln{\frac{1-T_2/T_1}{1+T_2/T_1}},
\label{symTHINC_d}
\end{equation}
where
\begin{equation*}
T_1=\tanh\left(\frac{\beta}{2}\right),\ T_2=\tanh\left(\frac{\breve{\alpha}_i\beta}{2}\right),\ \breve{\alpha}_i=\frac{\bar{q}_{i}-{q}_{a}}{{q}_{d}}.
\end{equation*}
Eventually, we get the symmetry-preserving formulae to compute the cell boundary values as
\begin{eqnarray}
\begin{aligned}
&\breve{q}_{i+\frac{1}{2}}^{L,T}=\breve{\mathcal{Q}}_i^{T}\left(x_{i+\frac{1}{2}}\right)={q}_{a}+{q}_{d}\frac{T_1+T_2/T_1}{1+T_2}, \\
&\breve{q}_{i-\frac{1}{2}}^{R,T}=\breve{\mathcal{Q}}_i^{T}\left(x_{i-\frac{1}{2}}\right)={q}_{a}-{q}_{d}\frac{T_1-T_2/T_1}{1-T_2}.
\label{symTHINC_boundary}
\end{aligned}
\end{eqnarray}
\begin{remark}
It is straightforwardly provable that \eqref{symTHINC_boundary} satisfy conditions   \eqref{cond_flip} and \eqref{cond_invers}, and thus ensures the symmetry in reconstruction.
\end{remark}

\begin{remark}
As shown in Fig. \ref{fig:THINC_scaling}, the reference point of the local coordinate $X_i$ (represented in blue) is shifted to the center of  cell $\Omega_i$, and the constant term in the THINC function (represented as the baseline in orange) is changed from $\bar{q}_{\mathrm{min}}$ to the average of values of the neighboring cells $(\bar{q}_{i-1}+\bar{q}_{i+1})/2$. These modifications enforce condition \eqref{cond_flip}. 
\end{remark}

\begin{remark}
The formula to obtain $\breve{\alpha}_i$ is different from the original $\alpha_i$ in Eq. \eqref{THINC_alpha} for symmetry property, where the small number $\epsilon$ is removed. Actually, $\epsilon$ in Eq. \eqref{THINC_alpha} is not necessary for avoiding zero-division which is already ruled out by condition $(\bar{q}_i-\bar{q}_{i-1})(\bar{q}_{i+1}-\bar{q}_i)>10^{-20}$. It can be straightforwardly shown that inclusion of  $\epsilon$ will violate condition \eqref{cond_invers}, and thus break the symmetry in numerical solutions.  

It is also noted that $\breve{\alpha}_i$ can be written in another form as
\begin{equation}
\breve{\alpha}_i=-\frac{\bar{q}_{i-1}-2\bar{q}_{i}+\bar{q}_{i+1}}{\bar{q}_{i+1}-\bar{q}_{i-1}}, 
\label{symTHINC_alpha2}
\end{equation}
which can be viewed as an approximation to the ratio of the second derivative and the first derivative of physical quantity $q$. However, this form of $\breve{\alpha}_i$ might cause the symmetry error due to the lack of associativity of summation of the numerator in \eqref{symTHINC_alpha2},  i.e. $(\bar{q}_{i-1}-2\bar{q}_{i})+\bar{q}_{i+1}\neq(\bar{q}_{i+1}-2\bar{q}_{i})+\bar{q}_{i-1}$.
\end{remark}

\begin{remark}
A more conventionally used alternative to \eqref{symTHINC_boundary} is to calculate $\breve{d}_i$ from \eqref{symTHINC_d} first, and then substitute it into the reconstruction function \eqref{symTHINC_function} to obtain the boundary values. However, this option does not ensure  \eqref{cond_flip} and \eqref{cond_invers}, as  for a logarithm function the following equality does not always hold  
\begin{eqnarray}
\ln{(a)}=-\ln{\left(\frac{1}{a}\right)}. 
\end{eqnarray}
In general, numerical processing of particular functions does  follow the mathematical commutativity.  
\end{remark}
It should be noted that \eqref{symTHINC_boundary} is not the only scheme that preserves symmetry in the THINC reconstruction though, it is a reasonable choice with simplicity and efficiency.

\subsection{HLLC Riemann solver}
\label{subsec:HLLC}

Another important component is the approximate Riemann solver. In this subsection, we examine symmetry-breaking causes of the HLLC Riemann solver \cite{Toro1994,Toro2009} used in the   P$_4$T$_2$-BVD method for Euler equations.

HLLC Riemann solver is the modification of HLL Riemann solver \cite{Harten1983a} to restore the contact wave. The numerical flux of the HLLC Riemann solver is expressed as
\begin{align}
\mathbf{F}^{HLLC}=\left\{
\begin{alignedat}{2}
&\mathbf{F}^L,& &\text{if}\ \ 0\leq s^L,\\
&\mathbf{F}^{*L},& &\text{if}\ \ s^L\leq 0\leq s^*,\\
&\mathbf{F}^{*R},&\quad &\text{if}\ \ s^*\leq 0\leq s^R,\\
&\mathbf{F}^R,& &\text{if}\ \ s^R\leq 0,
\end{alignedat}
\right.
\label{HLLCflux_choice}
\end{align}
where the asterisk ($*$) means the intermediate region between the left and right waves of  speeds $s^L$ and $s^R$ respectively. The middle contact wave moves at speed $s^*$, and divides the intermediate region further into two parts. From the Rankine-Hugoniot condition, the flux function and the conservative variable in the intermediate region can be obtained by 
\begin{align}
&\mathbf{F}^{*K}=\mathbf{F}^K+s^K(\mathbf{U}^{*K}-\mathbf{U}^K), \\
&\mathbf{U}^{*K}=\frac{s^K-u^K}{s^K-s^*}
\begin{pmatrix}
\rho^K \\ \rho^Ks^* \\ E^K+(s^*-u^K)\left(\rho^Ks^*+\frac{p^K}{s^K-u^K}\right)
\end{pmatrix}
,
\end{align}
with $K=L\ \text{or}\ R$. The intermediate wave speed $s^*$ is calculated as
\begin{align}
s^*=\frac{p^R-p^L+\rho^Lu^L(s^L-u^L)-\rho^Ru^R(s^R-u^R)}{\rho^L(s^L-u^L)-\rho^R(s^R-u^R)}.
\label{HLLC_sstar}
\end{align}
Here, the wave speeds $s^L$ and $s^R$ are obtained by the pressure-based wave speed estimation \cite{Toro1994}, 
\begin{align}
&s^L=u^L-c^Lq^L,\\ &s^R=u^R+c^Rq^R,
\end{align}
where
\begin{align}
q^K=\left\{
\begin{alignedat}{2}
&1,& &\text{if}\ \ p^*\leq p^K,\\
&\sqrt{1+\frac{\gamma +1}{2\gamma}\left( \frac{p^*}{p^K}-1\right)},&\quad &\text{otherwise},
\end{alignedat}
\right.
\end{align}
with $K=L\ \text{or}\ R$. The pressure in intermediate region is approximated by PVRS method \cite{Toro2009} as follows,
\begin{align}
p^*=\text{max}\left( 0,\ \frac{1}{2}(p^L+p^R)-\frac{1}{2}(u^R-u^L)\bar{\rho}\bar{c}\right),
\end{align}
where
\begin{align}
\bar{\rho}=\frac{1}{2}(\rho^L+\rho^R),\ \bar{c}=\frac{1}{2}(c^L+c^R).
\end{align}

There are three symmetry-breaking spots in the numerical formulation of HLLC Riemann solver. 

The first is in the flux function for transporting the transverse momentum component, which may break the diagonal symmetry. When HLLC solver is implemented dimension-wisely, the following formulae of flux functions may cause symmetry error,
\begin{align}
\mathbf{F}=
\begin{pmatrix}
\rho u \\ \rho u^2+p \\ \boxed{\rho uv} \\ (E+p)u \\
\end{pmatrix}
,\ \mathbf{G}=
\begin{pmatrix}
\rho v \\ \boxed{\rho uv} \\ \rho v^2+p \\ (E+p)v \\
\end{pmatrix}, 
\label{flux-hllc}
\end{align}
where the terms to advect the transverse momentums are boxed.  

Because the values of $u$ and $v$ are interchanged in dimension-wise calculations, \eqref{flux-hllc} might break the diagonal symmetry relationship due to the lack of the associativity of multiplication in numerical processing shown in Eq. \eqref{lack_of_associativity_multiple}. More specifically, the third component of $\mathbf{F}$ does not match the second component of $\mathbf{G}$ as the multiplications are arranged in different sequences, 
\begin{align}
{F}_3=(\underset{\text{a}}{\uwave{\rho}}\times\underset{\text{b}}{\uwave{u}})\times\underset{\text{c}}{\uwave{v}} \neq
(\underset{\text{a}}{\uwave{\rho}}\times\underset{\text{c}}{\uwave{u}})\times\underset{\text{b}}{\uwave{v}}={G}_2
\label{f3g2}
\end{align}
where the variable labeled by the same letter needs to be placed in the same order. Hence, the order of the multiplication in \eqref{f3g2} should be matched as Eq. \eqref{flux_function} or use the following form,
\begin{align}
\mathbf{\breve{F}}=
\begin{pmatrix}
\rho u \\ \rho u^2+p \\ \boxed{\rho vu} \\ (E+p)u \\
\end{pmatrix}
,\ \mathbf{\breve{G}}=
\begin{pmatrix}
\rho v \\ \boxed{\rho uv} \\ \rho v^2+p \\ (E+p)v \\
\end{pmatrix}
.
\label{HLLCflux_component}
\end{align}

The second point concerns the choice of the final numerical flux in Eq. \eqref{HLLCflux_choice}. When $s^*=0$, choosing either $\mathbf{F}^{*L}$ or $\mathbf{F}^{*R}$ will break the symmetry because they are not numerically identical due to the rounding error. This problem can be solved by adding the case of $s^*=0$ separately as follows,
\begin{align}
\breve{\mathbf{F}}^{HLLC}=\left\{
\begin{alignedat}{2}
&\mathbf{F}^L,& &\text{if}\ \ 0\leq s^L,\\
&\mathbf{F}^{*L},& &\text{if}\ \ s^L\leq 0< s^*,\\
&(\mathbf{F}^{*L}+\mathbf{F}^{*R})/2,&\quad &\text{if}\ \ s^*=0,\\
&\mathbf{F}^{*R},& &\text{if}\ \ s^*< 0\leq s^R,\\
&\mathbf{F}^R,& &\text{if}\ \ s^R\leq 0. 
\end{alignedat}
\right.
\label{HLLCflux_choice_sym}
\end{align}
It is noted that the following formulation proposed in \cite{Johnsen2006,Fleischmann2019} is effective to resolve this problem.
\begin{align}
\breve{\mathbf{F}}^{HLLC}=\frac{1+\mathrm{sgn}(s^*)}{2}\left\{\mathbf{F}^L+s^-\left(\mathbf{U}^{*L}-\mathbf{U}^L\right)\right\}+\frac{1-\mathrm{sgn}(s^*)}{2}\left\{\mathbf{F}^R+s^+\left(\mathbf{U}^{*R}-\mathbf{U}^R\right)\right\},
\label{HLLCflux_choice_sym_compact}
\end{align}
where $s^-=\mathrm{min}(s^L,0)$ and $s^+=\mathrm{max}(s^R,0)$.

The third possible cause of symmetry-breaking is the formula to calculate the speed of the intermediate wave $s^*$ shown in \eqref{HLLC_sstar}. Since the numerator of $s^*$ is composed of the summation of four terms, the order of the summation affects the symmetry property. A symmetry-preserving version of $s^*$ has been proposed in \cite{Fleischmann2019} as 
\begin{align}
\breve{s}^*=\frac{p^R-p^L+\pmb{\Bigl(}\rho^Lu^L(s^L-u^L)-\rho^Ru^R(s^R-u^R)\pmb{\Bigl)}}{\rho^L(s^L-u^L)-\rho^R(s^R-u^R)}.
\label{HLLC_sstar_sym}
\end{align}
The modification is adding a bracket to the original formula.

Our numerical experiments show that the symmetry errors in HLLC Riemann solver can be eliminated by modifying the formulations to fix the three possible causes analyzed above.

Before end this section, we summarize the symmetry-preserving techniques for the P$_4$T$_2$-BVD finite volume solver as follows.  
\begin{enumerate}[i)]

\item \label{item_W_to_U_con} Change the order of summation in the transformation from the characteristic variables to the conservative variables, using \eqref{solution_char_y-axis} or \eqref{eigen_sym}; 

\item \label{item_U_to_W_con} Add a bracket between the terms of $\bar{\rho u}$ and $\bar{\rho v}$ in the transformation from conservative variables to the characteristic variables using \eqref{solution_char_diag_x} and \eqref{solution_char_diag_y}; 

\item \label{item_reconstruction_con} Replace the formulae to compute the left-side and right-side values at cell boundaries in the P$_4$T$_2$-BVD scheme with the following symmetry-preserving formulae:
\begin{enumerate}[\ref{item_reconstruction_con}-i)]
    \item for the 4th-degree polynomial function use  \eqref{P4_boundary_sym}  or \eqref{P4_boundary_sym_lr},
    \item for the THINC function use \eqref{symTHINC_boundary};
\end{enumerate}
\item \label{item_HLLC_con} Modify the formulae in HLLC Riemann solver that break symmetry property:
\begin{enumerate}[\ref{item_HLLC_con}-i)]
    \item \label{item_HLLC_component} rearrange the order of multiplication in the advection flux for the transverse momentum component as in \eqref{HLLCflux_component},
    \item \label{item_HLLC_choice} use Eq. \eqref{HLLCflux_choice_sym_compact} to calculate numerical fluxes,  
    \item \label{item_HLLC_sstar} use  Eq. \eqref{HLLC_sstar_sym} to compute the speed of intermediate wave $s^*$.     
\end{enumerate}
\end{enumerate}

\section{Numerical results}
\label{sec:Results}

To verify the symmetry-preserving techniques introduced in section \ref{sec:Symmetry}, we have simulated some benchmark tests that have symmetric solution structures in space and are sensitive to any possible causes for symmetry breaking. The results demonstrate that the proposed techniques are effective for preserving both axis and diagonal symmetry properties. The 2D inviscid Euler equations are solved with relatively high mesh resolutions so that the asymmetric disturbances caused by floating-point arithmetic may grow and eventually contaminate the numerical solution.

The reconstruction is conducted by the P$_4$T$_2$-BVD scheme in terms of the characteristic variables. Symmetrized HLLC Riemann solver explained in subsection \ref{subsec:HLLC} is used to compute the numerical fluxes. Time evolution is performed by third-order SSP Runge-Kutta method \cite{Gottlieb2005}. The CFL number is set to 0.6 for all benchmark tests. The calculations have been done on Intel Xeon CPU E5-2687W 0 @ 3.10 GHz in a multi-threaded fashion on CentOS 6.10 operating system. GCC compiler version 4.4.6 was used without any special optimization options in compile commands.

\subsection{2D Riemann problems}

In 2D Riemann problems, the computational domain is divided into four areas and each area is given constant values of physical quantities initially. In some cases,  Kelvin-Helmholtz instabilities are produced along the contact discontinuities dividing the domain. The development of the vortices is usually used as the measure to evaluate the numerical dissipation of a scheme. We choose configurations 3 and 12 from the 2D Riemann problem test set \cite{Kurganov2002}, which have diagonal symmetry in solution structures.

\begin{figure}[htbp]
    \centering
    \includegraphics[width=8.0cm]{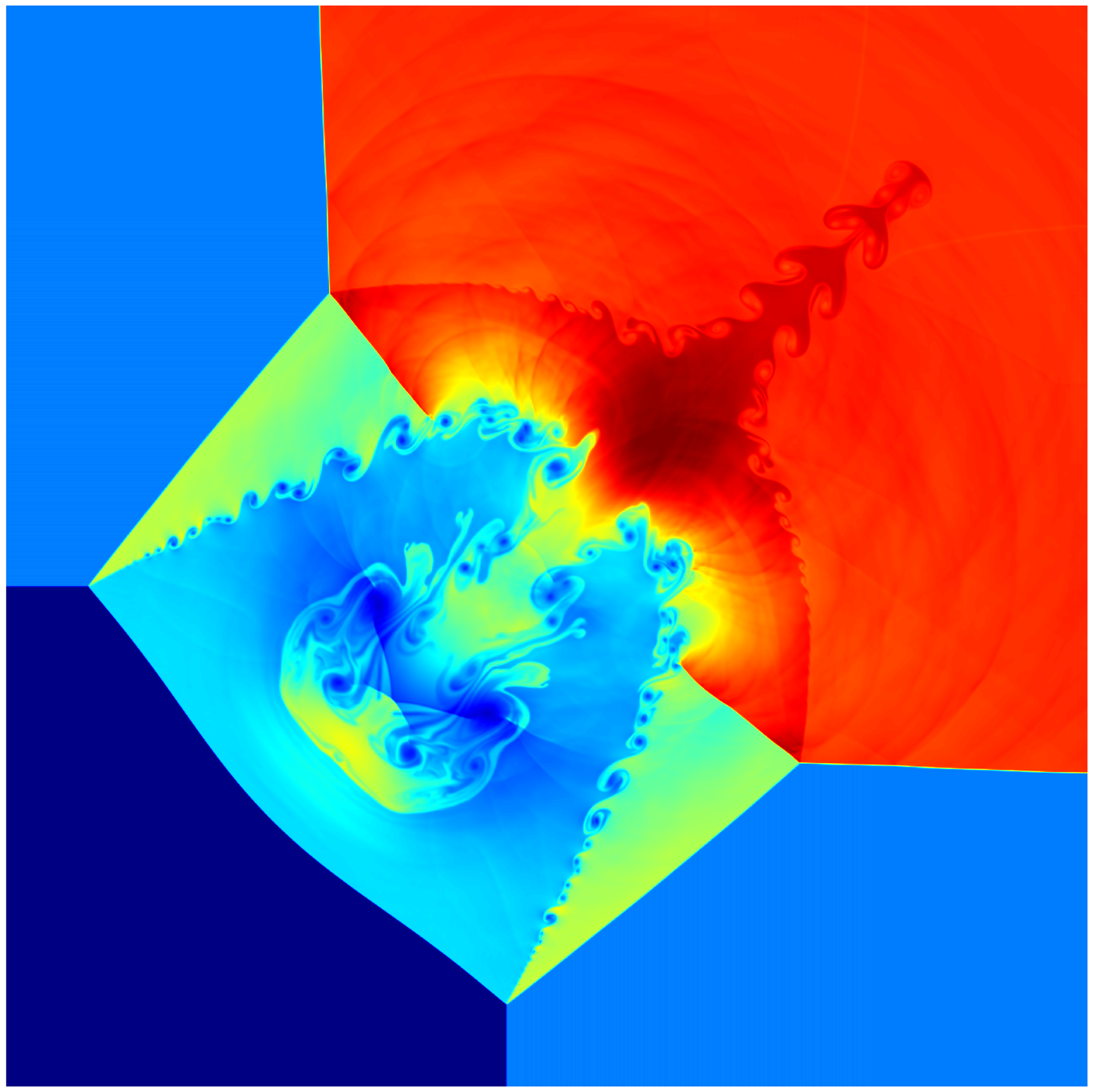}
    \includegraphics[width=8.0cm]{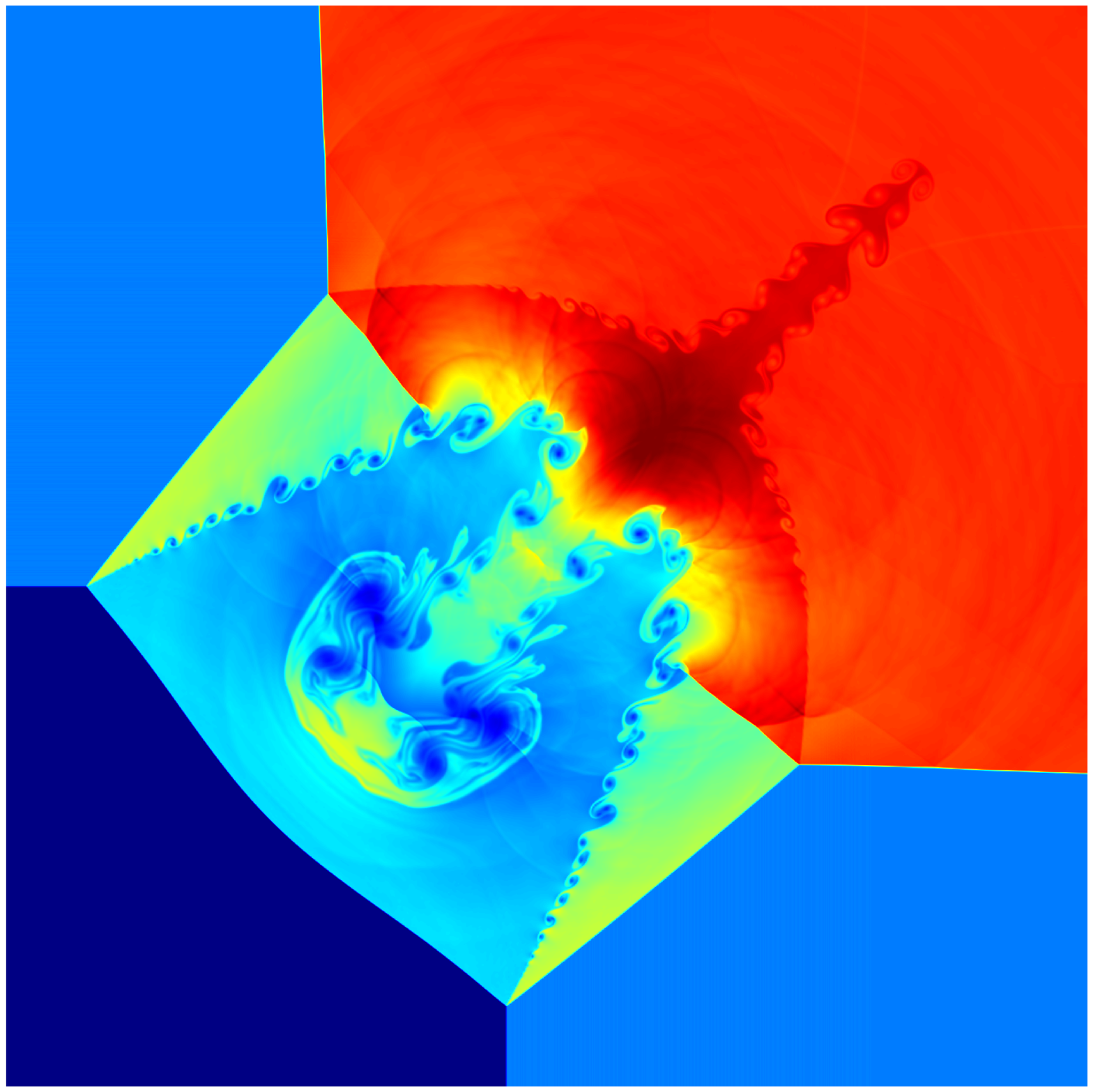}
    \caption{Numerical solutions of density for config. 3 of 2D Riemann problems (blue=0.135 to red=1.75). Left panel is for the original scheme and right panel is for the symmetry-preserving scheme.}
    \label{fig:Riemann_3}
\end{figure}

\begin{figure}[htbp]
    \centering
    \includegraphics[width=8.0cm]{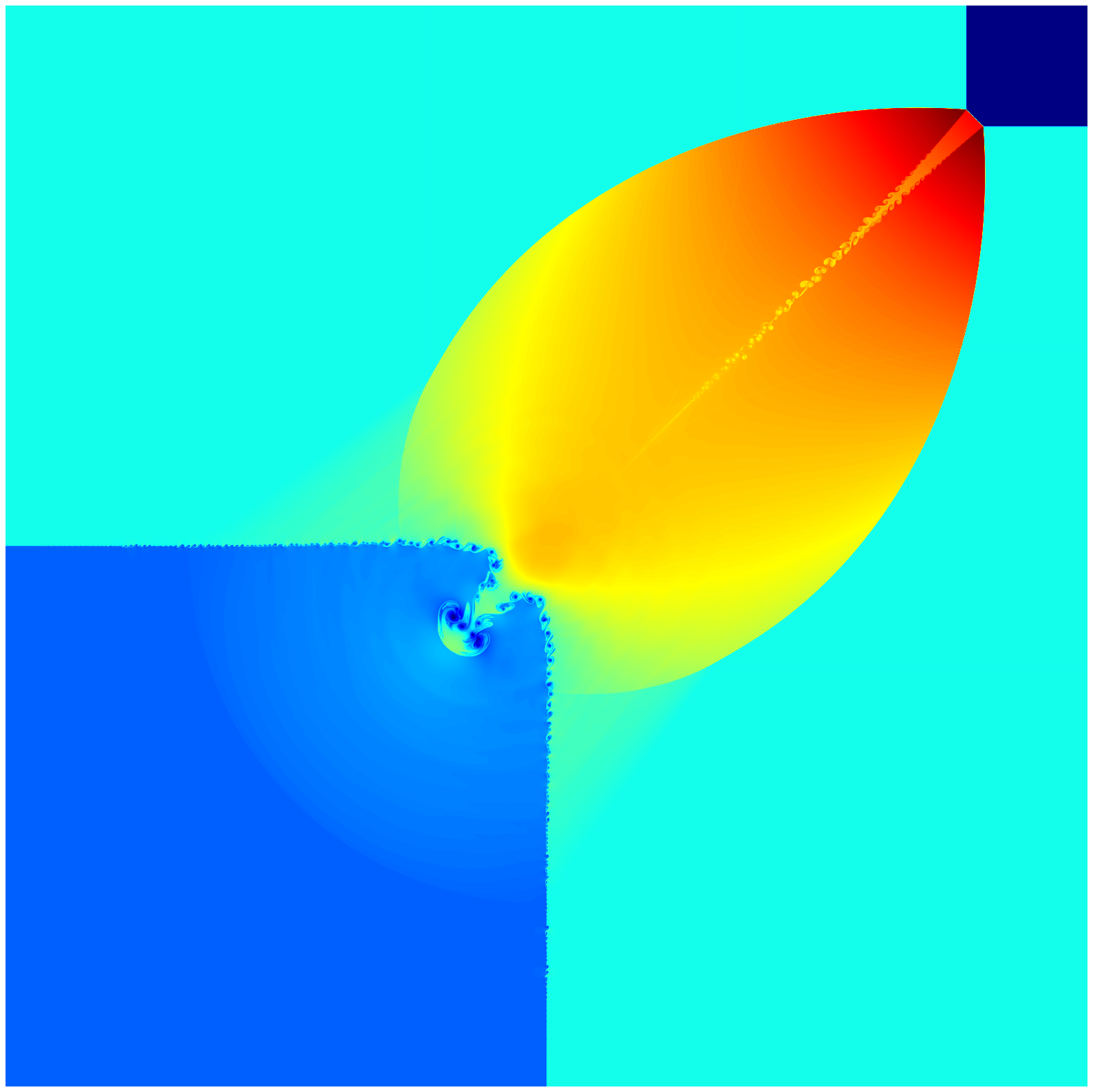}
    \includegraphics[width=8.0cm]{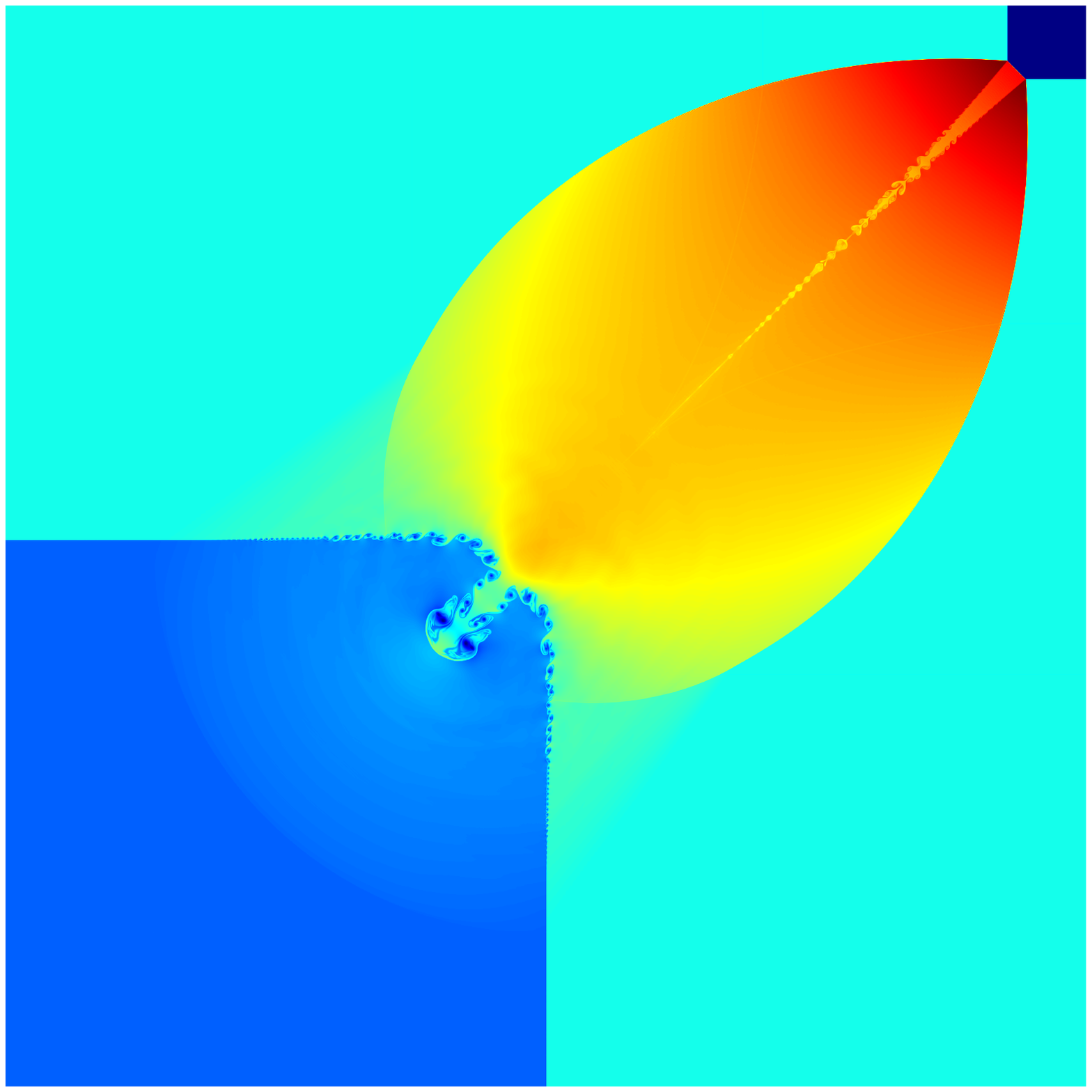}
    \includegraphics[width=8.0cm]{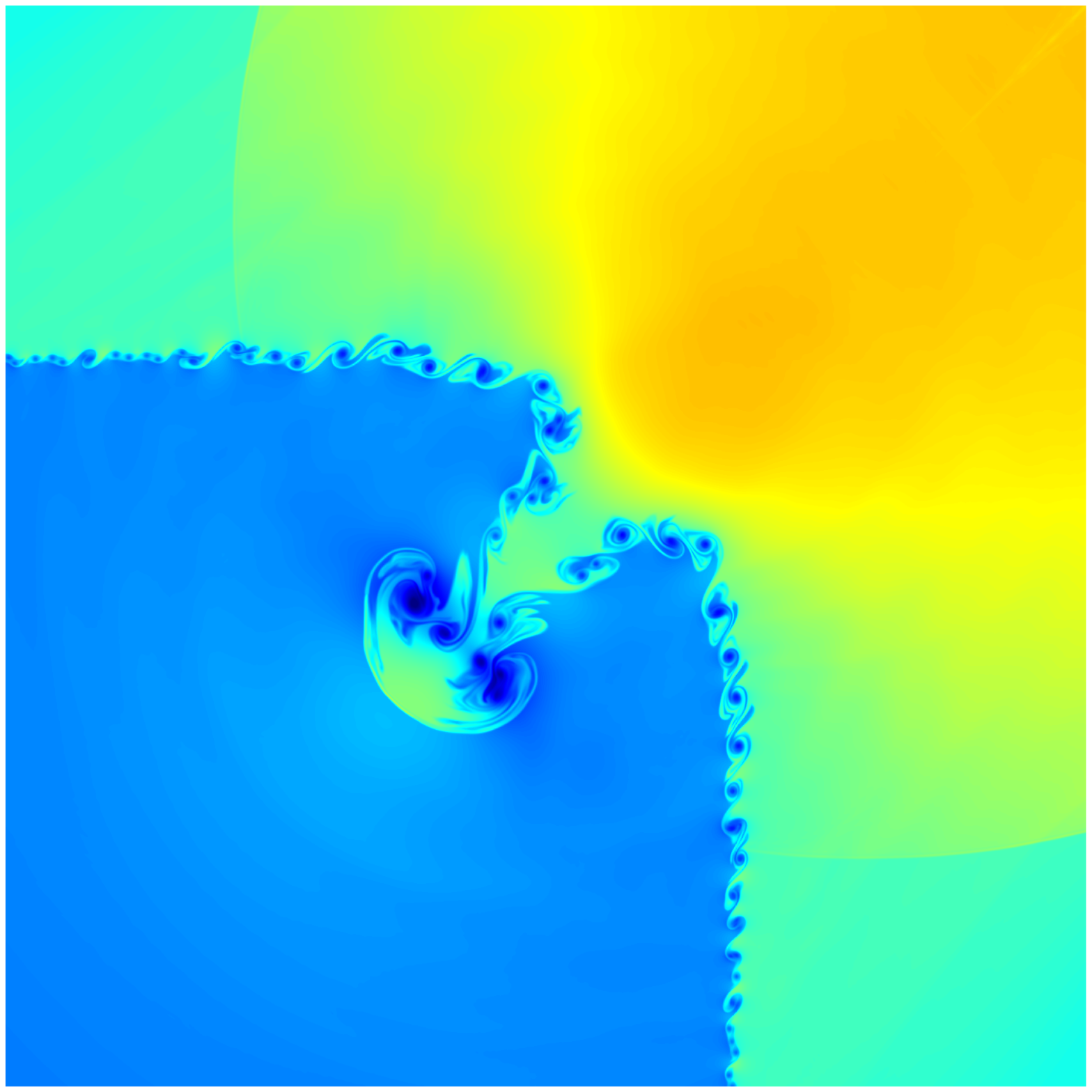}
    \includegraphics[width=8.0cm]{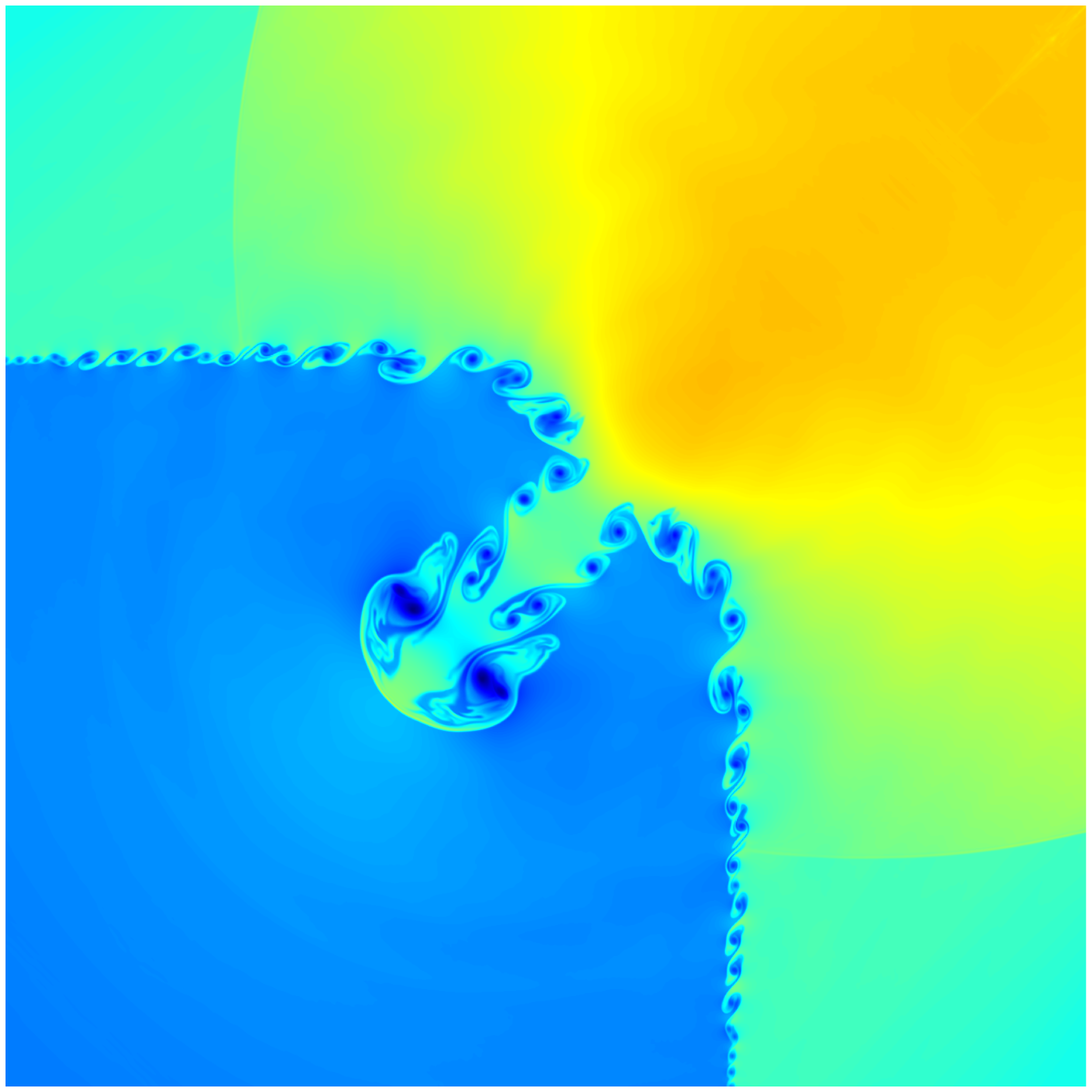}
    \caption{Numerical solutions of density for config. 12 of 2D Riemann problems (blue=0.55 to red=1.7). Left panels show the results from the original scheme and right panels are for the symmetry-preserving scheme. The bottom figures show the zoomed region of $[-0.2,0.1]\times[-0.2,0.1]$.}
    \label{fig:Riemann_12}
\end{figure}

The initial conditions of config. 3 and 12 are as follows,
\begin{align}
&\text{Config.\ 3:} \ \ (\rho_0,u_0,v_0,p_0)=\left\{
\begin{alignedat}{3}
&(1.5,0,0,1.5),&\quad \ &x>0.3-\varepsilon,&\ &y>0.3-\varepsilon, \\
&(0.5323,1.206,0,0.3),& &x<0.3-\varepsilon,& &y>0.3+\varepsilon, \\
&(0.138,1.206,1.206,0.029),& &x<0.3+\varepsilon,& &y<0.3+\varepsilon, \\
&(0.5323,0,1.206,0.3),& &x>0.3+\varepsilon,& &y<0.3-\varepsilon,
\end{alignedat}
\right.
\label{ini_Riemann_3}
\\
&\text{Config.\ 12:}\ \ (\rho_0,u_0,v_0,p_0)=\left\{
\begin{alignedat}{3}
&(0.5313,0,0,0.4),&\quad \ &x>-\varepsilon,&\ &y>-\varepsilon, \\
&(1,0.7276,0,1),& &x<-\varepsilon,& &y>\varepsilon, \\
&(0.8,0,0,1),& &x<\varepsilon,& &y<\varepsilon, \\
&(1,0,0.7276,1),& &x>\varepsilon,& &y<-\varepsilon,
\end{alignedat}
\right.
\label{ini_Riemann_12}
\end{align}
where the small number $\varepsilon=10^{-15}$ is introduced to avoid symmetry errors in the initial conditions. The boundary condition is set to zero-gradient for all boundaries. The computational domain is $[-0.5,0.5]\times[-0.5,0.5]$ and partitioned by a $1000\times 1000$ grid for config. 3 and a $4000\times 4000$ grid for config. 12. The numerical results of the density at time $t=0.8$ of config. 3 are shown in Fig. \ref{fig:Riemann_3} and at time $t=0.25$ of config. 12 are shown in Fig. \ref{fig:Riemann_12}.

It is observed from the left panels in Figs. \ref{fig:Riemann_3} and \ref{fig:Riemann_12} that the results of the original P$_4$T$_2$-BVD scheme cannot preserve the diagonal symmetry property, especially in the vortex structures along slip lines. On the contrary, the symmetry property is perfectly preserved by the symmetry-preserving techniques introduced in this paper, as shown in the right panels of Figs. \ref{fig:Riemann_3} and \ref{fig:Riemann_12}. These results indicate that the proposed techniques are valid to preserve the diagonal symmetry.

\subsection{Rayleigh-Taylor instability}

The Rayleigh-Taylor instability (RTI) is a typical benchmark test that has axis-symmetricity in solution structure. In the past, RTI has been extensively used to assess numerical methods for Euler equations \cite{Xu2005, Don2018,Fu2018}, and most existing schemes result in asymmetric solutions in high-resolution simulations, where lower numerical dissipation allows asymmetric rounding errors to develop and be visible in numerical solution. 

The initial condition is set by placing a heavy fluid above a light fluid under gravity with a velocity perturbation on the interface between the two fluids, as specified below,
\begin{align}
(\rho_0,u_0,v_0,p_0)=
\begin{cases}
(2,0,v_0(x),2y+1), & \text{if}\ \ y<0.5, \\
(1,0,v_0(x),y+1.5), & \text{otherwise},
\end{cases}
\label{ini_Rayleigh}
\end{align}
where $v_0(x)$ is the initial value of the velocity in $y$-direction perturbed in $x$-direction to trigger the mixing between heavy and light fluids. Conventionally, $v_0(x)$ is given by 
\begin{align}
v_0(x)=-0.025c\cos{(8\pi x)}. 
\label{ini_Rayleigh_asym}
\end{align}
However, as pointed out by Fleischmann et al. \cite{Fleischmann2019}, $v_0(x)$ in \eqref{ini_Rayleigh_asym} is not symmetric due to the fact that $\cos{(\pi-\epsilon)}\neq\cos{(\pi+\epsilon)}$ in the floating-point arithmetic. Fleischmann et al. \cite{Fleischmann2019} proposed the following symmetry-preserving formulation for specifying $v_0(x)$, 
\begin{align}
\breve{v}_0(x)=
\begin{cases}
-0.025c\cos{(8\pi x)}, & \text{if}\ \ x<0.125, \\
-0.025c\cos{(0.25-8\pi x)}, & \text{otherwise}.
\end{cases}
\label{ini_Rayleigh_sym}
\end{align}

\begin{figure}[htbp]
    \centering
    \includegraphics[width=4.0cm]{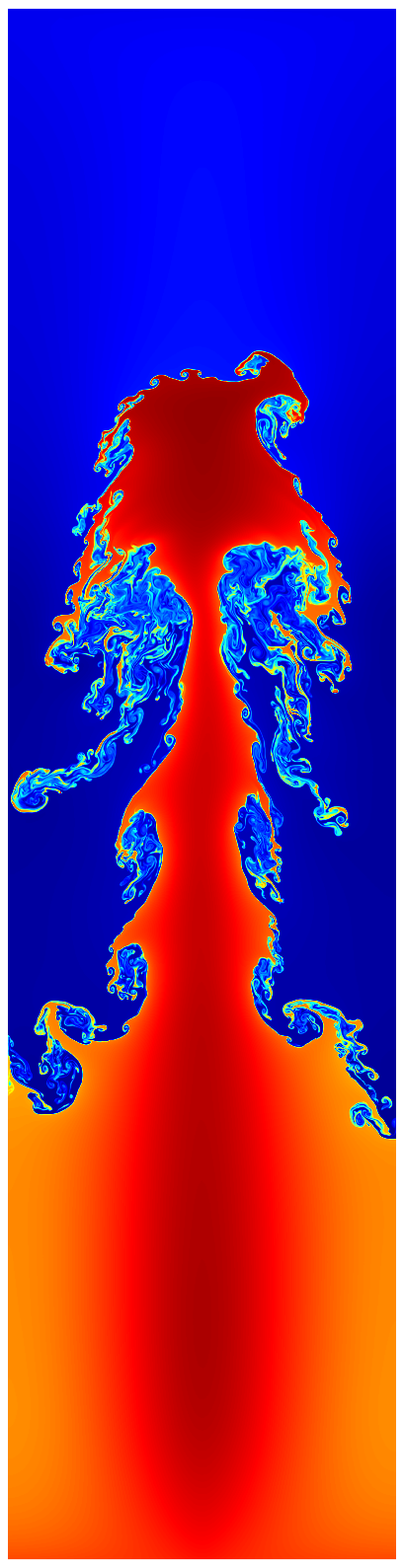}
    \includegraphics[width=4.0cm]{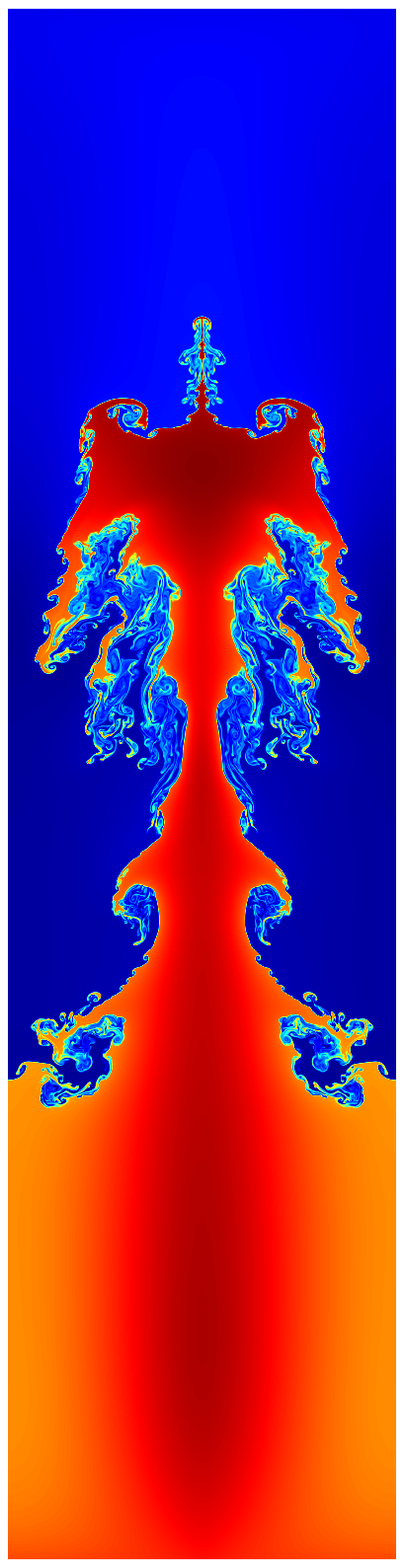}
    \caption{Numerical solutions of density for Rayleigh-Taylor instability (blue=0.85 to red=2.25) with the mesh resolution of $1024\times 4096$. Left panel is the result of the original scheme and right panel is that of the symmetry-preserving scheme.}
    \label{fig:Rayleigh_1024}
\end{figure}

\begin{figure}[htbp]
    \centering
    \begin{minipage}[b]{4cm}
    \includegraphics[width=4.0cm]{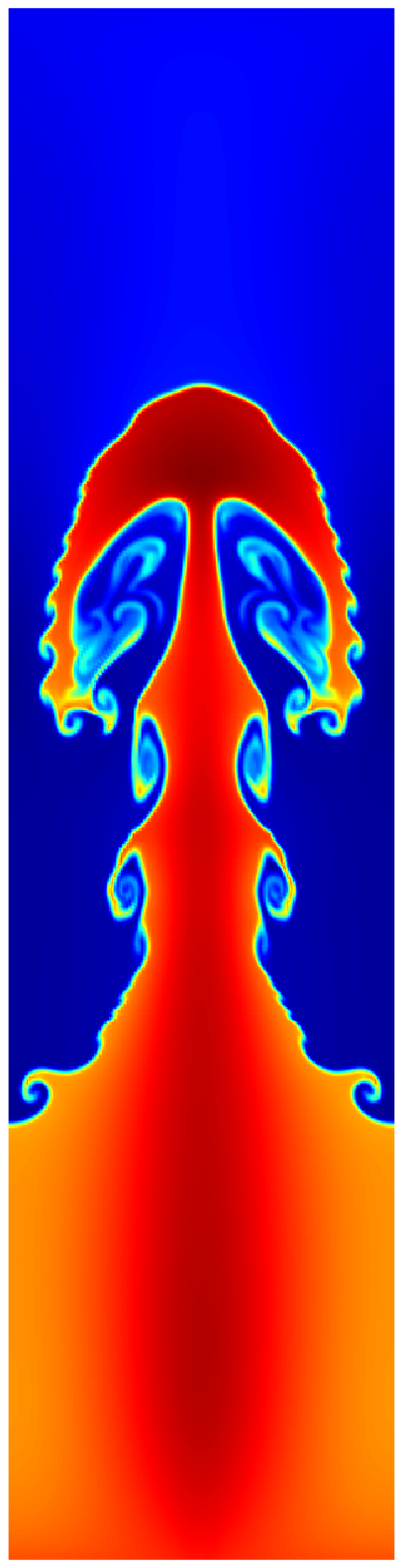}
    \end{minipage}
    \begin{minipage}[b]{8.4cm}
    \includegraphics[width=2.0cm]{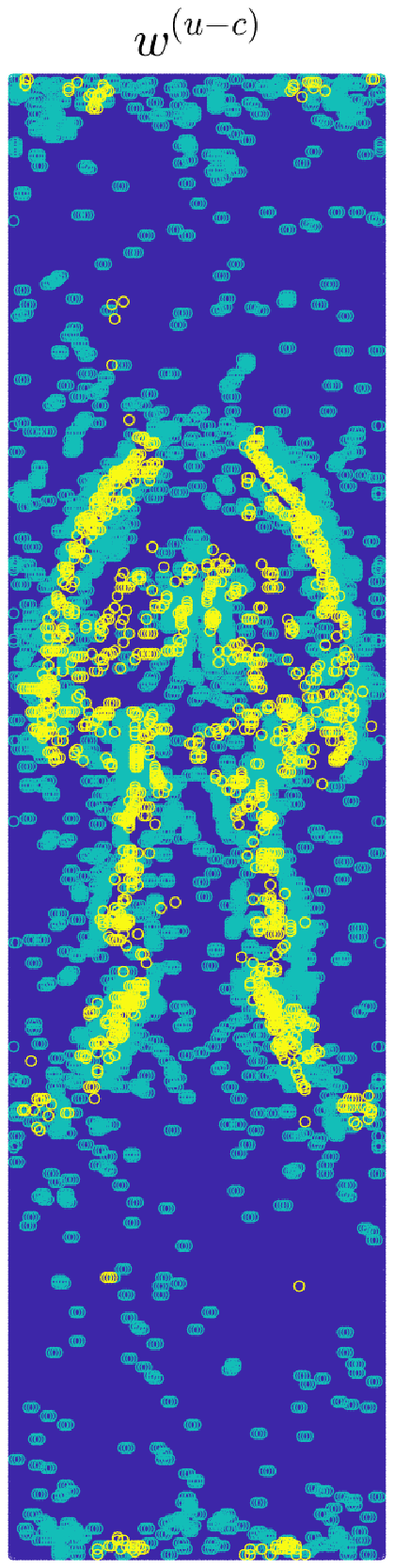}
    \includegraphics[width=2.0cm]{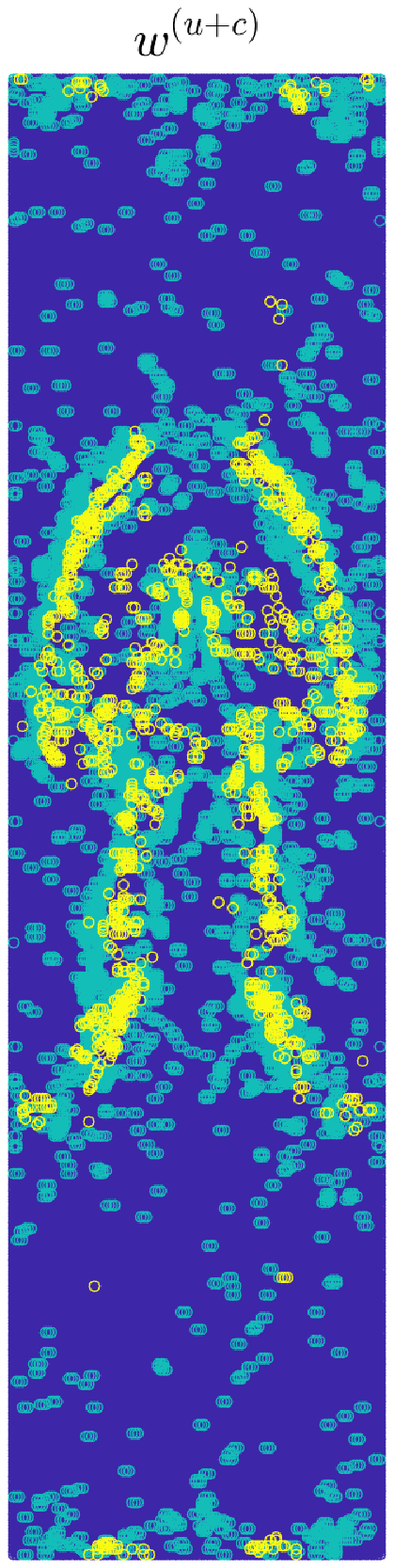}
    \includegraphics[width=2.0cm]{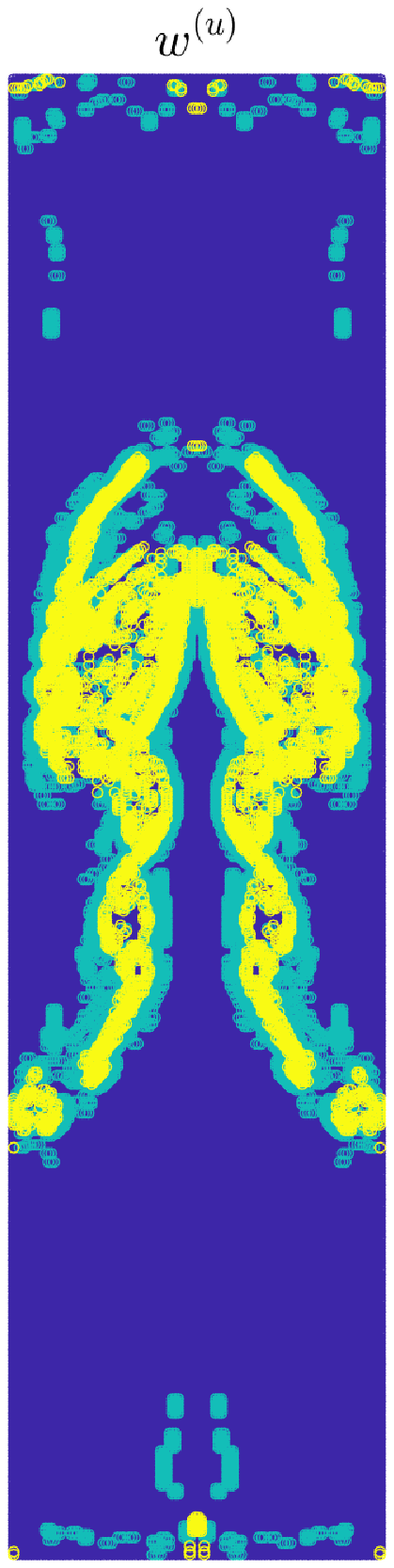}
    \includegraphics[width=2.0cm]{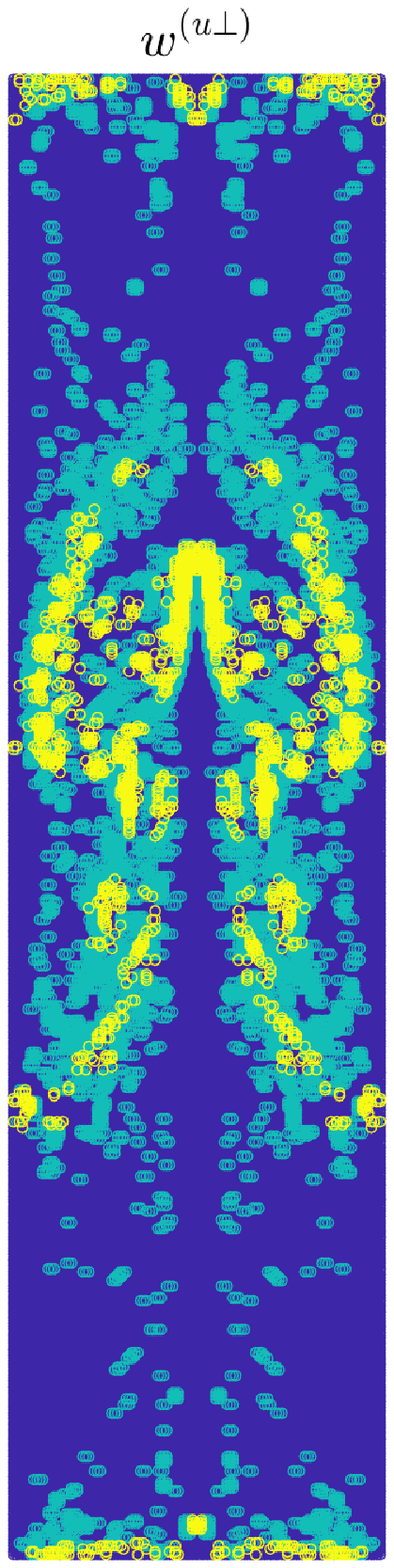}
    \includegraphics[width=2.0cm]{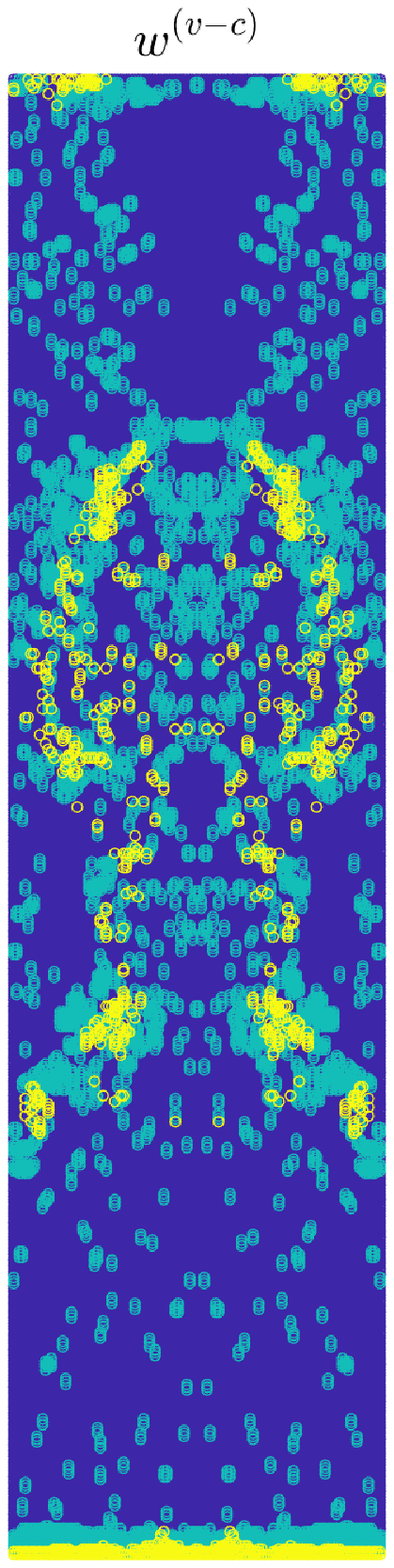}
    \includegraphics[width=2.0cm]{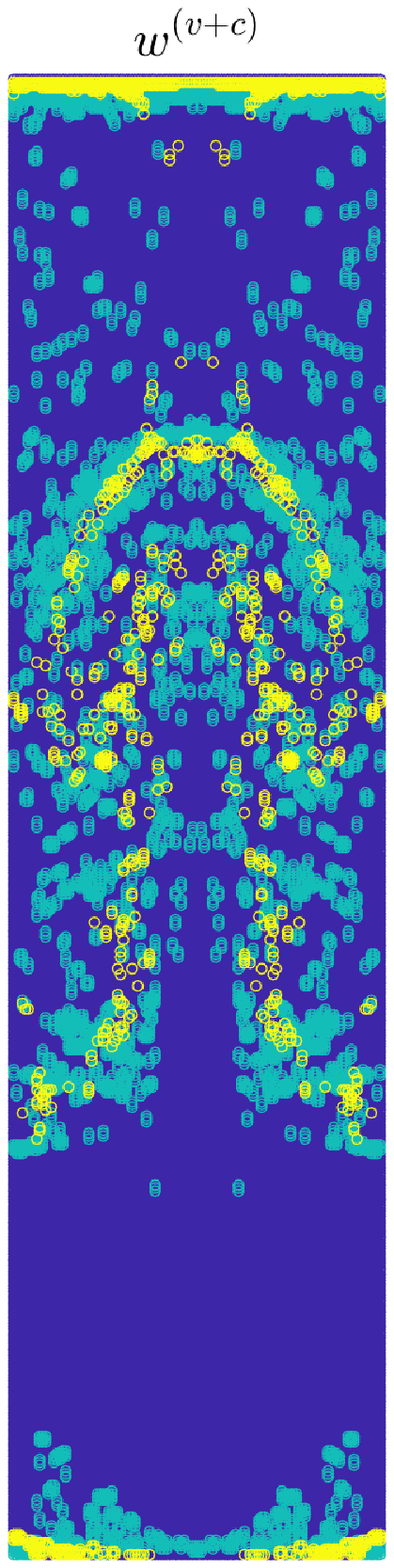}
    \includegraphics[width=2.0cm]{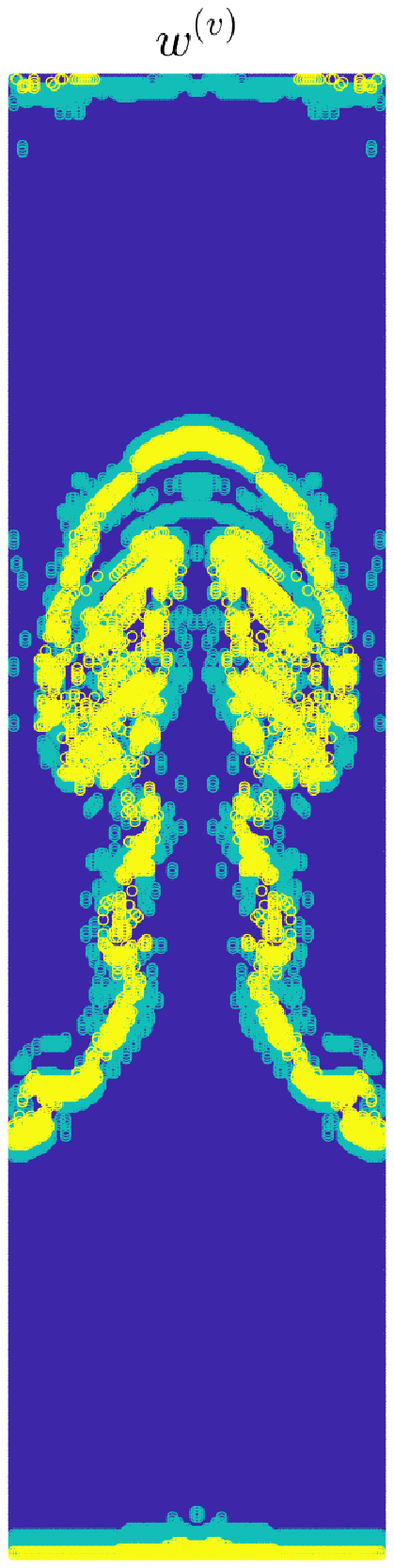}
    \includegraphics[width=2.0cm]{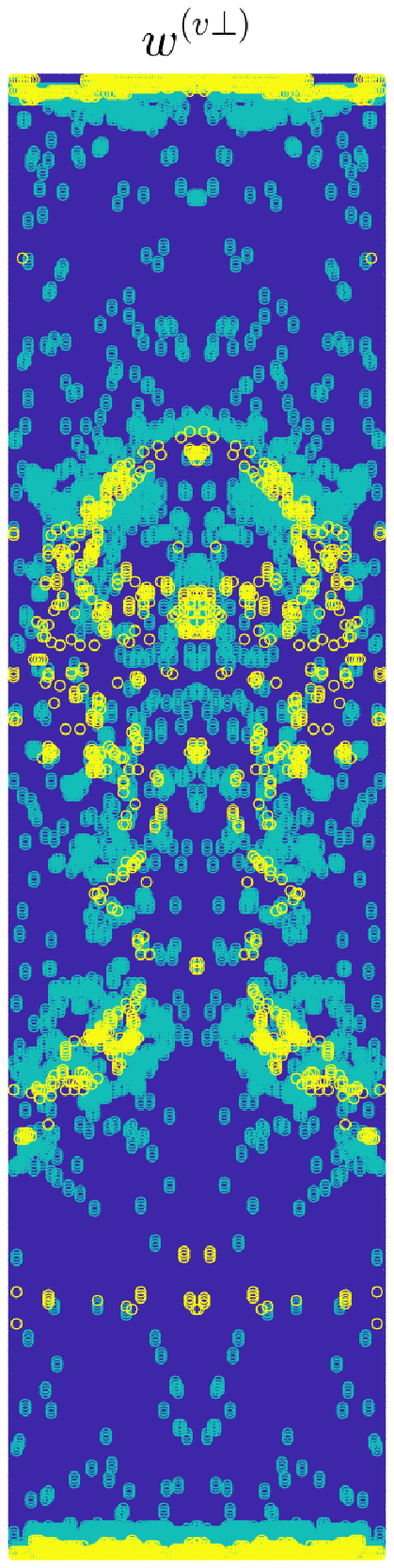}
    \end{minipage}
    \caption{Left panel is the numerical solutions of density for Rayleigh-Taylor instability calculated by symmetry-preserving scheme (blue=0.85 and red=2.25) and right panel shows which reconstruction function is selected by the BVD algorithm for each characteristic variable (purple=4th-degree polynomial, light blue=THINC($\beta_s$), yellow=THINC($\beta_l$)). Mesh resolution is $128\times 512$.}
    \label{fig:Rayleigh_128_distW}
\end{figure}

Reflective boundary condition is imposed at the left and right boundaries and the top and bottom boundaries are fixed as $(\rho,u,v,p)_{top}=(1,0,0,2.5)$ and $(\rho,u,v,p)_{bottom}=(2,0,0,1)$ respectively. The computational domain is $[0,0.25]\times[0,1]$. The specific heat ratio is $\gamma=\frac{5}{3}$. The numerical results of the density at time $t=1.95$ are shown in Figs. \ref{fig:Rayleigh_1024} and \ref{fig:Rayleigh_128_distW}.

Fig. \ref{fig:Rayleigh_1024} shows that the $y$-axis asymmetric flow structure calculated by the original scheme is improved, and perfect symmetry is reproduced by the symmetry-preserving techniques even in this low-dissipation simulation on a high-resolution mesh. This result indicates that the proposed methods are effective for preserving axis symmetry. 

We further examine the symmetric property in the BVD algorithm to select the reconstruction function. Fig. \ref{fig:Rayleigh_128_distW} shows which candidate function is selected in the BVD algorithm as the final interpolant for reconstruction. Mesh cells are marked with different colors according to which reconstruction function is used. I.e. purple indicates a cell where the 4th-degree polynomial is used, light blue for a cell using   THINC($\beta_s$), and yellow for THINC($\beta_l$). It can be seen that THINC($\beta_l$) is selected in the cells in the vicinity of strong discontinuous solutions while THINC($\beta_s$) is chosen for the cells where the solution is between smooth and discontinuous. The selection of reconstruction function in the BVD algorithm is completely symmetric for all characteristic variables except $w^{(u-c)}$ and $w^{(u+c)}$. As shown in Fig. \ref{fig:dist_Wx} and table \ref{tab:relationships_summary},  $w^{(u-c)}$ and $w^{(u+c)}$ have their values interchanged with symmetric positions in the $y$-axis symmetry.  Consistently, Figs. \ref{fig:Rayleigh_1024} and \ref{fig:Rayleigh_128_distW} shows that the left and right halves of the selected function of $w^{(u-c)}$ are swapped with the right and left halves of $w^{(u+c)}$ inversely. It completely agrees with the $y$-axis symmetry of the characteristic variables analyzed in section \ref{sec:Symmetry}.

\subsection{Implosion test}

This test demonstrates an implosion phenomenon where a diamond-shaped low-pressure region is crushed by the surrounding high-pressure \cite{Liska2003}. The implosion test has been used to verify the performance of the numerical codes for astrophysical simulation \cite{Sutherland2010,Schneider2015}. Since the diagonal jets are sensitive to numerical disturbances, the solution tends to easily deviate from the diagonal direction if the code contains any symmetry error. It is a very challenging benchmark test to evaluate both axis and diagonal symmetry-preserving properties of numerical methods.  

\begin{figure}[htbp]
    \centering
    \includegraphics[width=8.0cm]{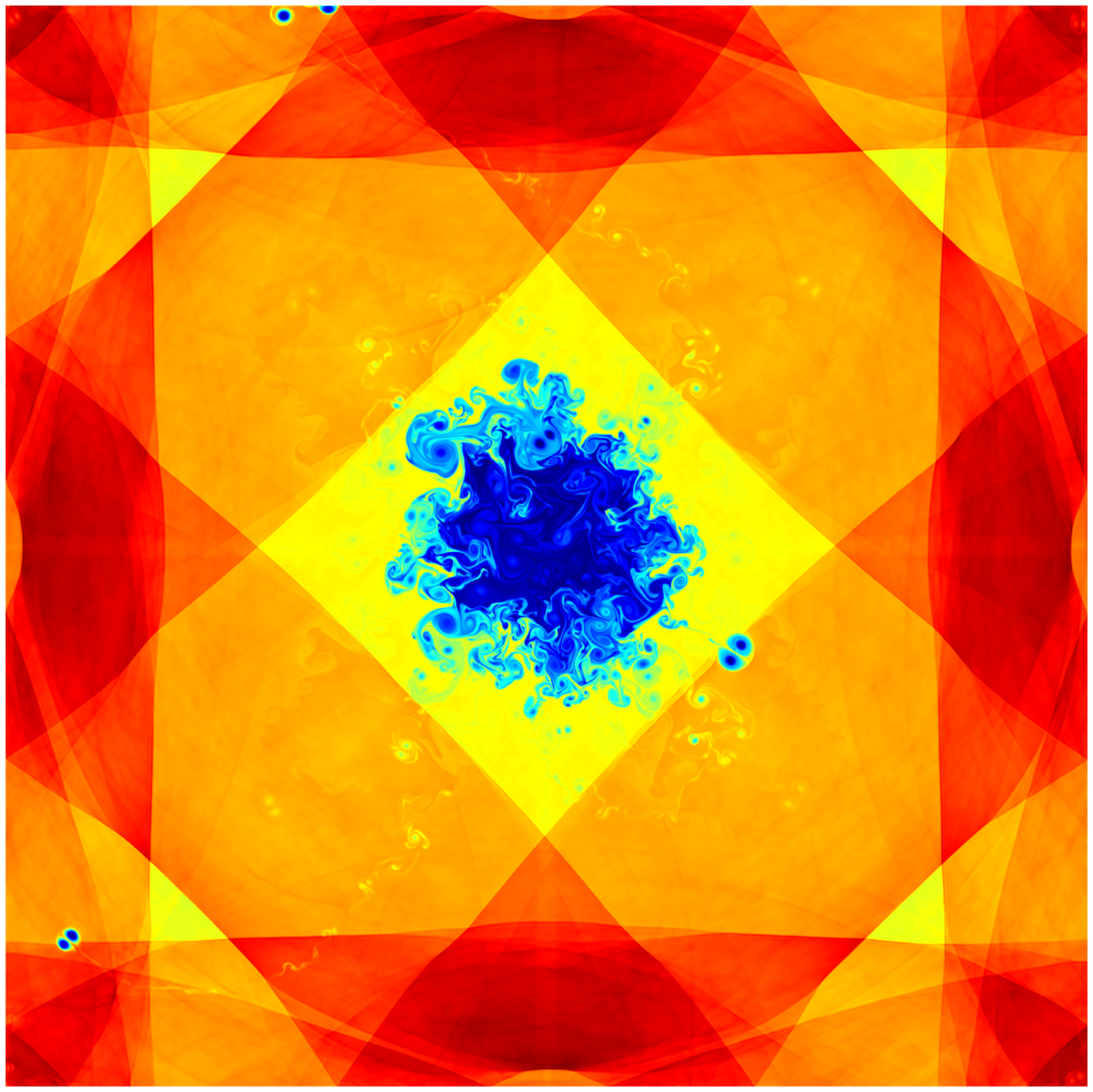}
    \includegraphics[width=8.0cm]{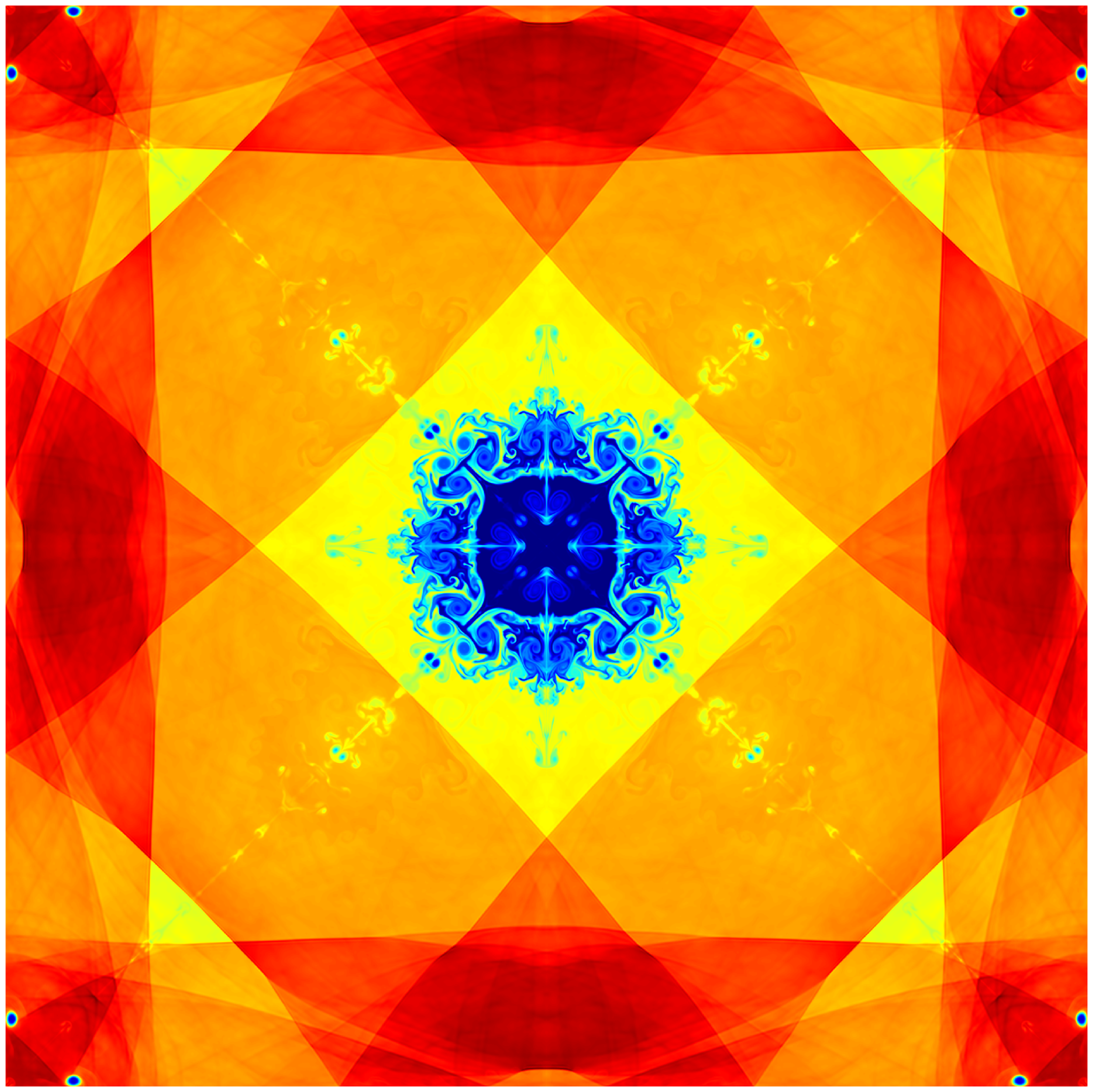}
    \includegraphics[width=8.0cm]{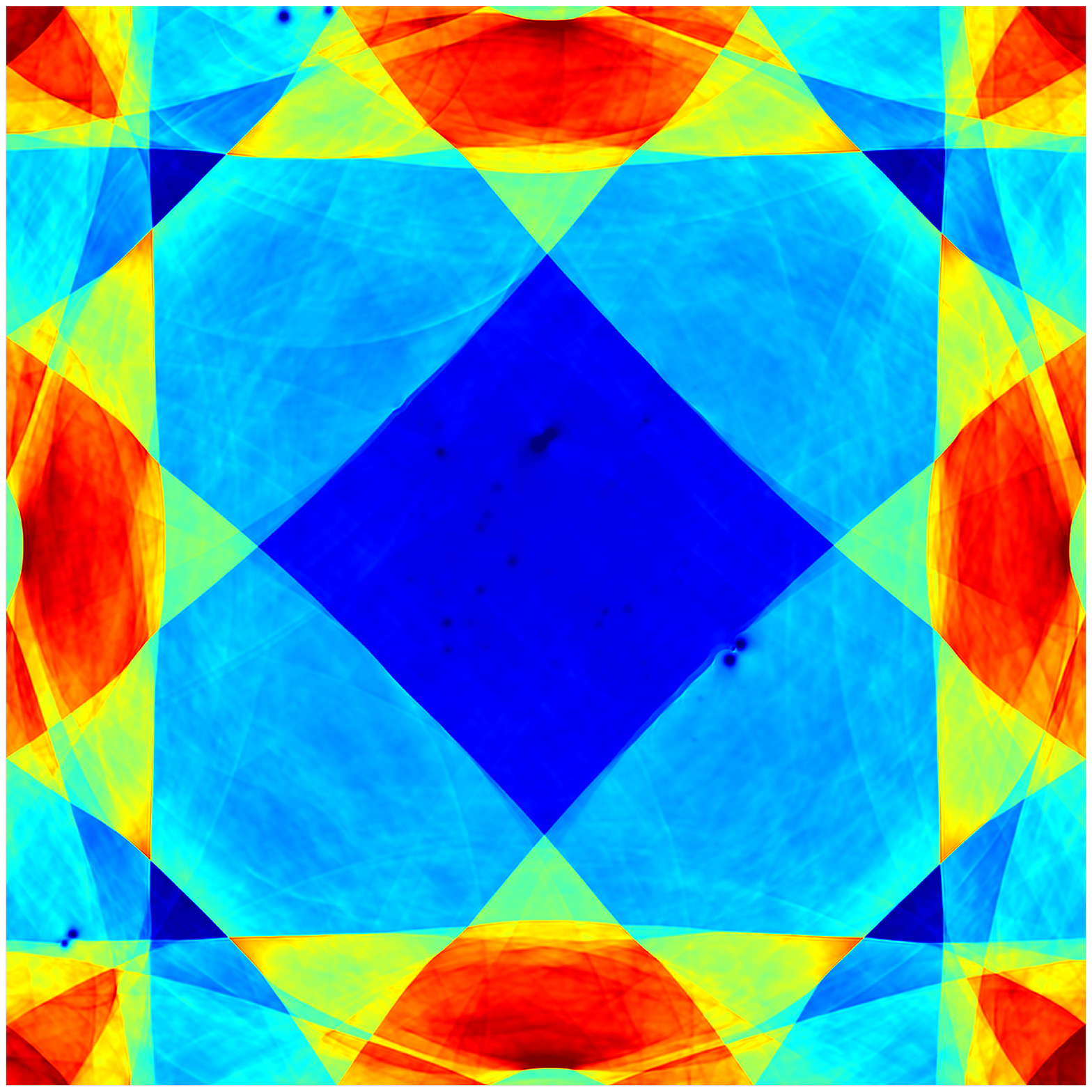}
    \includegraphics[width=8.0cm]{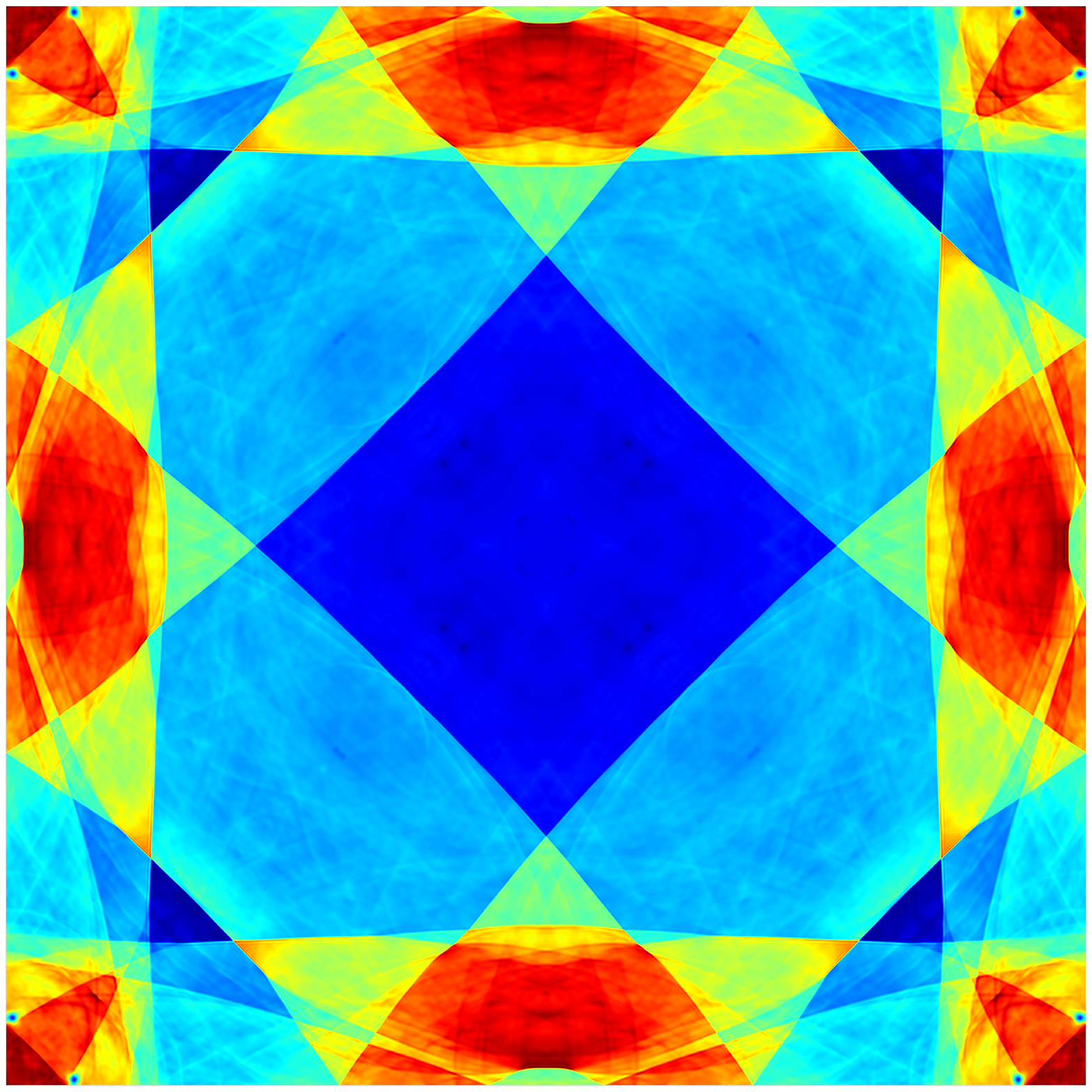}
    \caption{Numerical solutions for implosion test with resolution $1600\times 1600$. Left panels are the results of the original scheme and right panels are for the symmetry-preserving scheme. The density (blue=0.41 to red=1.07) is on the top, and the pressure (blue=0.74 to red=1.1) is on the bottom.}
    \label{fig:implosion_1600}
\end{figure}

The initial conditions are set as follows,
\begin{align}
(\rho_0,u_0,v_0,p_0)=\left\{
\begin{alignedat}{2}
&(0.125,0,0,0.14),&\quad \ &|y+x|<0.15+\varepsilon,\ |y-x|<0.15+\varepsilon, \\
&(1,0,0,1),& &\text{otherwise},
\end{alignedat}
\right.
\label{ini_implosion}
\end{align}
where $\varepsilon=10^{-10}$ is introduced in order to avoid asymmetricity in initial condition. The reflective boundary condition is imposed on all boundaries of computational domain which is a square area specified by $[-0.3,0.3]\times[-0.3,0.3]$. The numerical results of the density and the pressure at time $t=2.5$ are shown in Fig. \ref{fig:implosion_1600}.

It can be seen that the trajectories of the jets in the results of the original scheme are distorted because of the symmetry errors. By adopting the proposed techniques, both axis and diagonal symmetric properties are perfectly preserved.

\section{Conclusion}
\label{sec:Conclusion}

In this study, the mechanisms of symmetry-breaking in a finite volume framework with P$_4$T$_2$-BVD reconstruction scheme are thoroughly examined. As the high-resolution simulations using high-order numerical schemes on fine computational grids might amplify the symmetry errors of round-off level which break the symmetric solution structures, attention has been paid to remove the possible causes due to the lack of associativity in floating-point arithmetic. 

The reasons for symmetry-breaking exist in different components of the numerical solver for Euler equations. They are mainly attributed to 1) the sequence of arithmetic operations in transformation calculations between conservative and characteristic variables, 2) the arrangement of the stencil for spatial reconstruction, 3) the commutative operations with nonlinear functions, and 4) the asymmetricity associated with the “if” logic operation. 

Numerical modifications and numerical techniques are proposed to completely remove the possible causes for symmetry breaking in  the P$_4$T$_2$-BVD finite volume solver. Benchmark tests which have symmetric solution structures are used to verify the methods proposed in this paper. The numerical results demonstrate that the numerical solver can perfectly reproduce the symmetric solution structures in all tests. 


\section*{Acknowledgment}

This work was supported in part by the fund from JSPS (Japan Society for the Promotion of Science) under Grant Nos. 18H01366 and 19H05613.

\appendix
\setcounter{table}{0}
\setcounter{figure}{0}
\setcounter{equation}{0}

\section{The diagonal symmetricity in the transformation from characteristic variables to conservative variables}

 The transformations of the reconstructed characteristic variables $\mathbf{W}$ to the conservative variables $\mathbf{U}$ at the diagonally symmetric positions are computed by $\mathbf{U}_A=\mathbf{R}_{xA}\cdot\mathbf{W}_{xA}$ and $\mathbf{U}_B=\mathbf{R}_{yB}\cdot\mathbf{W}_{yB}$ in $x$- and $y$-directions respectively. For point $A$, we have
\begin{align}
\begin{pmatrix}
\rho \\ \rho u \\ \rho v \\ E
\end{pmatrix}
_A&=
\begin{pmatrix}
1 \\ u-c \\ v \\ H-uc
\end{pmatrix}
_Aw^{(u-c)}_A+
\begin{pmatrix}
1 \\ u \\ v \\ \frac{u^2+v^2}{2}
\end{pmatrix}
_Aw^{(u)}_A+
\begin{pmatrix}
1 \\ u+c \\ v \\ H+uc
\end{pmatrix}
_Aw^{(u+c)}_A+
\begin{pmatrix}
0 \\ 0 \\ 1 \\ v
\end{pmatrix}
_Aw^{(u\perp)}_A. 
\end{align}

Using the diagonal symmetry rule given in \eqref{prim_diag} and table \ref{tab:relationships_summary}, we get the transformation formulae at point $B$ as follows, 
\begin{align}
\begin{pmatrix}
\rho \\ \rho u \\ \rho v \\ E
\end{pmatrix}
_B&=
\begin{pmatrix}
1 \\ u \\ v-c \\ H-vc
\end{pmatrix}
_Bw^{(v-c)}_B+
\begin{pmatrix}
1 \\ u \\ v \\ \frac{u^2+v^2}{2}
\end{pmatrix}
_Bw^{(v)}_B+
\begin{pmatrix}
1 \\ u \\ v+c \\ H+vc
\end{pmatrix}
_Bw^{(v+c)}_B+
\begin{pmatrix}
0 \\ 1 \\ 0 \\ u
\end{pmatrix}
_Bw^{(v\perp)}_B \nonumber \\ &=
\begin{pmatrix}
1 \\ v \\ u-c \\ H-uc 
\end{pmatrix}
_Aw^{(u-c)}_A+
\begin{pmatrix}
1 \\ v \\ u \\ \frac{v^2+u^2}{2}
\end{pmatrix}
_Aw^{(u)}_A+
\begin{pmatrix}
1 \\ v \\ u+c \\ H+uc
\end{pmatrix}
_Aw^{(u+c)}_A+
\begin{pmatrix}
0 \\ 1 \\ 0 \\ v
\end{pmatrix}
_Aw^{(u\perp)}_A \nonumber \\ &=
\begin{pmatrix}
\rho \\ \rho v \\ \rho u \\ E
\end{pmatrix}
_A.
\label{WB_to_UB_diag}
\end{align}
It is observed that all terms are computed in the same order when applying the transformation dimension-wisely in $x$- and $y$-directions. So, no any modification at this stage is required to enhance the symmetry in numerical solution. 

\bibliography{Bibliography}

\begin{thebibliography}{59}
\expandafter\ifx\csname natexlab\endcsname\relax\def\natexlab#1{#1}\fi
\providecommand{\url}[1]{\texttt{#1}}
\providecommand{\href}[2]{#2}
\providecommand{\path}[1]{#1}
\providecommand{\DOIprefix}{doi:}
\providecommand{\ArXivprefix}{arXiv:}
\providecommand{\URLprefix}{URL: }
\providecommand{\Pubmedprefix}{pmid:}
\providecommand{\doi}[1]{\href{http://dx.doi.org/#1}{\path{#1}}}
\providecommand{\Pubmed}[1]{\href{pmid:#1}{\path{#1}}}
\providecommand{\bibinfo}[2]{#2}
\ifx\xfnm\relax \def\xfnm[#1]{\unskip,\space#1}\fi
\bibitem[{Godunov(1959)}]{Godunov1959}
\bibinfo{author}{S.~K. Godunov},
\newblock \bibinfo{title}{{A difference method for numerical calculation of
  discontinuous solutions of the equations of hydrodynamics}},
\newblock \bibinfo{journal}{Math. Sb.} \bibinfo{volume}{47}
  (\bibinfo{year}{1959}) \bibinfo{pages}{271--306}.
\bibitem[{Harten(1983)}]{Harten1983}
\bibinfo{author}{A.~Harten},
\newblock \bibinfo{title}{{High Resolution Schemes for Hyperbolic Conservation
  Laws}},
\newblock \bibinfo{journal}{International Journal for Numerical Methods in
  Fluids} \bibinfo{volume}{23} (\bibinfo{year}{1983})
  \bibinfo{pages}{309--323}.
\bibitem[{{Van Leer}(1977)}]{VanLeer1977}
\bibinfo{author}{B.~{Van Leer}},
\newblock \bibinfo{title}{{Towards the ultimate conservative difference scheme.
  IV. A new approach to numerical convection}},
\newblock \bibinfo{journal}{Journal of Computational Physics}
  \bibinfo{volume}{23} (\bibinfo{year}{1977}) \bibinfo{pages}{276--299}.
\bibitem[{Harten et~al.(1987)Harten, Engquist, Osher, and
  Chakravarthy}]{Harten1987}
\bibinfo{author}{A.~Harten}, \bibinfo{author}{B.~Engquist},
  \bibinfo{author}{S.~Osher}, \bibinfo{author}{S.~R. Chakravarthy},
\newblock \bibinfo{title}{{Uniformly High Order Accurte Essentially
  Non-Oscillatory Schemes III}},
\newblock \bibinfo{journal}{Journal of Computational Physics}
  \bibinfo{volume}{71} (\bibinfo{year}{1987}) \bibinfo{pages}{231--323}.
\bibitem[{Shu and Osher(1988)}]{Shu1988}
\bibinfo{author}{C.-W. Shu}, \bibinfo{author}{S.~Osher},
\newblock \bibinfo{title}{{Efficient Implementation of Essentially
  Non-oscillatory Shock-Capturing Schemes}},
\newblock \bibinfo{journal}{Journal of Computational Physics}
  \bibinfo{volume}{77} (\bibinfo{year}{1988}) \bibinfo{pages}{439--471}.
\bibitem[{Liu et~al.(1994)Liu, Osher, and Chan}]{Liu1994}
\bibinfo{author}{X.-D. Liu}, \bibinfo{author}{S.~Osher},
  \bibinfo{author}{T.~Chan},
\newblock \bibinfo{title}{{Weighted Essentially Non-oscillatory Schemes}},
\newblock \bibinfo{journal}{Journal of Computational Physics}
  \bibinfo{volume}{115} (\bibinfo{year}{1994}) \bibinfo{pages}{200--212}.
\bibitem[{Jiang and Shu(1996)}]{Jiang1996}
\bibinfo{author}{G.-S. Jiang}, \bibinfo{author}{C.-W. Shu},
\newblock \bibinfo{title}{{Efficient implementation of weighted ENO schemes}},
\newblock \bibinfo{journal}{Journal of Computational Physics}
  \bibinfo{volume}{126} (\bibinfo{year}{1996}) \bibinfo{pages}{202--228}.
\bibitem[{Henrick et~al.(2005)Henrick, Aslam, and Powers}]{Henrick2005}
\bibinfo{author}{A.~K. Henrick}, \bibinfo{author}{T.~D. Aslam},
  \bibinfo{author}{J.~M. Powers},
\newblock \bibinfo{title}{{Mapped weighted essentially non-oscillatory schemes:
  Achieving optimal order near critical points}},
\newblock \bibinfo{journal}{Journal of Computational Physics}
  \bibinfo{volume}{207} (\bibinfo{year}{2005}) \bibinfo{pages}{542--567}.
\bibitem[{Borges et~al.(2008)Borges, Carmona, Costa, and Don}]{Borges2008}
\bibinfo{author}{R.~Borges}, \bibinfo{author}{M.~Carmona},
  \bibinfo{author}{B.~Costa}, \bibinfo{author}{W.~S. Don},
\newblock \bibinfo{title}{{An improved weighted essentially non-oscillatory
  scheme for hyperbolic conservation laws}},
\newblock \bibinfo{journal}{Journal of Computational Physics}
  \bibinfo{volume}{227} (\bibinfo{year}{2008}) \bibinfo{pages}{3191--3211}.
\bibitem[{Ha et~al.(2013)Ha, {Ho Kim}, {Ju Lee}, and Yoon}]{Ha2013}
\bibinfo{author}{Y.~Ha}, \bibinfo{author}{C.~{Ho Kim}}, \bibinfo{author}{Y.~{Ju
  Lee}}, \bibinfo{author}{J.~Yoon},
\newblock \bibinfo{title}{{An improved weighted essentially non-oscillatory
  scheme with a new smoothness indicator}},
\newblock \bibinfo{journal}{Journal of Computational Physics}
  \bibinfo{volume}{232} (\bibinfo{year}{2013}) \bibinfo{pages}{68--86}.
\bibitem[{Fan et~al.(2014)Fan, Shen, Tian, and Yang}]{Fan2014}
\bibinfo{author}{P.~Fan}, \bibinfo{author}{Y.~Shen}, \bibinfo{author}{B.~Tian},
  \bibinfo{author}{C.~Yang},
\newblock \bibinfo{title}{{A new smoothness indicator for improving the
  weighted essentially non-oscillatory scheme}},
\newblock \bibinfo{journal}{Journal of Computational Physics}
  \bibinfo{volume}{269} (\bibinfo{year}{2014}) \bibinfo{pages}{329--354}.
\bibitem[{Kim et~al.(2016)Kim, Ha, and Yoon}]{Kim2016}
\bibinfo{author}{C.~H. Kim}, \bibinfo{author}{Y.~Ha},
  \bibinfo{author}{J.~Yoon},
\newblock \bibinfo{title}{{Modified Non-linear Weights for Fifth-Order Weighted
  Essentially Non-oscillatory Schemes}},
\newblock \bibinfo{journal}{Journal of Scientific Computing}
  \bibinfo{volume}{67} (\bibinfo{year}{2016}) \bibinfo{pages}{299--323}.
\bibitem[{Fu et~al.(2016)Fu, Hu, and Adams}]{Fu2016}
\bibinfo{author}{L.~Fu}, \bibinfo{author}{X.~Y. Hu}, \bibinfo{author}{N.~A.
  Adams},
\newblock \bibinfo{title}{{A family of high-order targeted ENO schemes for
  compressible-fluid simulations}},
\newblock \bibinfo{journal}{Journal of Computational Physics}
  \bibinfo{volume}{305} (\bibinfo{year}{2016}) \bibinfo{pages}{333--359}.
\bibitem[{Acker et~al.(2016)Acker, Borges, and Costa}]{Acker2016}
\bibinfo{author}{F.~Acker}, \bibinfo{author}{R.~B. Borges},
  \bibinfo{author}{B.~Costa},
\newblock \bibinfo{title}{{An improved WENO-Z scheme}},
\newblock \bibinfo{journal}{Journal of Computational Physics}
  \bibinfo{volume}{313} (\bibinfo{year}{2016}) \bibinfo{pages}{726--753}.
\bibitem[{Fu et~al.(2017)Fu, Hu, and Adams}]{Fu2017}
\bibinfo{author}{L.~Fu}, \bibinfo{author}{X.~Y. Hu}, \bibinfo{author}{N.~A.
  Adams},
\newblock \bibinfo{title}{{Targeted ENO schemes with tailored resolution
  property for hyperbolic conservation laws}},
\newblock \bibinfo{journal}{Journal of Computational Physics}
  \bibinfo{volume}{349} (\bibinfo{year}{2017}) \bibinfo{pages}{97--121}.
\bibitem[{Rathan and {Naga Raju}(2018)}]{Rathan2018}
\bibinfo{author}{S.~Rathan}, \bibinfo{author}{G.~{Naga Raju}},
\newblock \bibinfo{title}{{A modified fifth-order WENO scheme for hyperbolic
  conservation laws}},
\newblock \bibinfo{journal}{Computers and Mathematics with Applications}
  \bibinfo{volume}{75} (\bibinfo{year}{2018}) \bibinfo{pages}{1531--1549}.
\bibitem[{Fu et~al.(2018)Fu, Hu, and Adams}]{Fu2018}
\bibinfo{author}{L.~Fu}, \bibinfo{author}{X.~Y. Hu}, \bibinfo{author}{N.~A.
  Adams},
\newblock \bibinfo{title}{{A new class of adaptive high-order targeted ENO
  schemes for hyperbolic conservation laws}},
\newblock \bibinfo{journal}{Journal of Computational Physics}
  \bibinfo{volume}{374} (\bibinfo{year}{2018}) \bibinfo{pages}{724--751}.
\bibitem[{Sun et~al.(2016)Sun, Inaba, and Xiao}]{Sun2016}
\bibinfo{author}{Z.~Sun}, \bibinfo{author}{S.~Inaba},
  \bibinfo{author}{F.~Xiao},
\newblock \bibinfo{title}{{Boundary Variation Diminishing (BVD) reconstruction:
  A new approach to improve Godunov schemes}},
\newblock \bibinfo{journal}{Journal of Computational Physics}
  \bibinfo{volume}{322} (\bibinfo{year}{2016}) \bibinfo{pages}{309--325}.
\bibitem[{Deng et~al.(2017)Deng, Xie, and Xiao}]{Deng2017}
\bibinfo{author}{X.~Deng}, \bibinfo{author}{B.~Xie}, \bibinfo{author}{F.~Xiao},
\newblock \bibinfo{title}{{A finite volume multi-moment method with boundary
  variation diminishing principle for Euler equation on three-dimensional
  hybrid unstructured grids}},
\newblock \bibinfo{journal}{Computers and Fluids} \bibinfo{volume}{153}
  (\bibinfo{year}{2017}) \bibinfo{pages}{85--101}.
\bibitem[{Xie et~al.(2017)Xie, Deng, Sun, and Xiao}]{Xie2017a}
\bibinfo{author}{B.~Xie}, \bibinfo{author}{X.~Deng}, \bibinfo{author}{Z.~Sun},
  \bibinfo{author}{F.~Xiao},
\newblock \bibinfo{title}{{A hybrid pressure–density-based Mach uniform
  algorithm for 2D Euler equations on unstructured grids by using multi-moment
  finite volume method}},
\newblock \bibinfo{journal}{Journal of Computational Physics}
  \bibinfo{volume}{335} (\bibinfo{year}{2017}) \bibinfo{pages}{637--663}.
\bibitem[{Deng et~al.(2018{\natexlab{a}})Deng, Inaba, Xie, Shyue, and
  Xiao}]{Deng2018a}
\bibinfo{author}{X.~Deng}, \bibinfo{author}{S.~Inaba},
  \bibinfo{author}{B.~Xie}, \bibinfo{author}{K.~M. Shyue},
  \bibinfo{author}{F.~Xiao},
\newblock \bibinfo{title}{{High fidelity discontinuity-resolving reconstruction
  for compressible multiphase flows with moving interfaces}},
\newblock \bibinfo{journal}{Journal of Computational Physics}
  \bibinfo{volume}{371} (\bibinfo{year}{2018}{\natexlab{a}})
  \bibinfo{pages}{945--966}.
\bibitem[{Deng et~al.(2018{\natexlab{b}})Deng, Xie, Loub{\`{e}}re, Shimizu, and
  Xiao}]{Deng2018}
\bibinfo{author}{X.~Deng}, \bibinfo{author}{B.~Xie},
  \bibinfo{author}{R.~Loub{\`{e}}re}, \bibinfo{author}{Y.~Shimizu},
  \bibinfo{author}{F.~Xiao},
\newblock \bibinfo{title}{{Limiter-free discontinuity-capturing scheme for
  compressible gas dynamics with reactive fronts}},
\newblock \bibinfo{journal}{Computers and Fluids} \bibinfo{volume}{171}
  (\bibinfo{year}{2018}{\natexlab{b}}) \bibinfo{pages}{1--14}.
\bibitem[{Deng et~al.(2019)Deng, Shimizu, and Xiao}]{Deng2019}
\bibinfo{author}{X.~Deng}, \bibinfo{author}{Y.~Shimizu},
  \bibinfo{author}{F.~Xiao},
\newblock \bibinfo{title}{{A fifth-order shock capturing scheme with two-stage
  boundary variation diminishing algorithm}},
\newblock \bibinfo{journal}{Journal of Computational Physics}
  \bibinfo{volume}{386} (\bibinfo{year}{2019}) \bibinfo{pages}{323--349}.
\bibitem[{Tann et~al.(2019)Tann, Deng, Shimizu, Loub{\`{e}}re, and
  Xiao}]{Tann2019}
\bibinfo{author}{S.~Tann}, \bibinfo{author}{X.~Deng},
  \bibinfo{author}{Y.~Shimizu}, \bibinfo{author}{R.~Loub{\`{e}}re},
  \bibinfo{author}{F.~Xiao},
\newblock \bibinfo{title}{{Solution property preserving reconstruction for
  finite volume scheme: a boundary variation diminishing+multidimensional
  optimal order detection framework}},
\newblock \bibinfo{journal}{International Journal for Numerical Methods in
  Fluids} \bibinfo{volume}{92} (\bibinfo{year}{2019})
  \bibinfo{pages}{603--634}.
\bibitem[{Deng et~al.(2020)Deng, Shimizu, Xie, and Xiao}]{Deng2020}
\bibinfo{author}{X.~Deng}, \bibinfo{author}{Y.~Shimizu},
  \bibinfo{author}{B.~Xie}, \bibinfo{author}{F.~Xiao},
\newblock \bibinfo{title}{{Constructing higher order discontinuity-capturing
  schemes with upwind-biased interpolations and boundary variation diminishing
  algorithm}},
\newblock \bibinfo{journal}{Computers and Fluids} \bibinfo{volume}{200}
  (\bibinfo{year}{2020}) \bibinfo{pages}{104433}.
\bibitem[{Tann et~al.(2020)Tann, Deng, Loub{\`{e}}re, and Xiao}]{Tann2020}
\bibinfo{author}{S.~Tann}, \bibinfo{author}{X.~Deng},
  \bibinfo{author}{R.~Loub{\`{e}}re}, \bibinfo{author}{F.~Xiao},
\newblock \bibinfo{title}{{Solution property preserving reconstruction BVD+MOOD
  scheme for compressible euler equations with source terms and detonations}},
\newblock \bibinfo{journal}{Computers and Fluids} \bibinfo{volume}{206}
  (\bibinfo{year}{2020}) \bibinfo{pages}{104594}.
\bibitem[{Cheng et~al.(2021)Cheng, Deng, Xie, Jiang, and Xiao}]{Cheng2021}
\bibinfo{author}{L.~Cheng}, \bibinfo{author}{X.~Deng},
  \bibinfo{author}{B.~Xie}, \bibinfo{author}{Y.~Jiang},
  \bibinfo{author}{F.~Xiao},
\newblock \bibinfo{title}{{Low-dissipation BVD schemes for single and
  multi-phase compressible flows on unstructured grids}},
\newblock \bibinfo{journal}{Journal of Computational Physics}
  \bibinfo{volume}{428} (\bibinfo{year}{2021}) \bibinfo{pages}{110088}.
\bibitem[{Jiang et~al.(2021)Jiang, Deng, Xiao, Yan, Yu, and Lou}]{Jiang2021}
\bibinfo{author}{Z.-H. Jiang}, \bibinfo{author}{X.~Deng},
  \bibinfo{author}{F.~Xiao}, \bibinfo{author}{C.~Yan}, \bibinfo{author}{J.~Yu},
  \bibinfo{author}{S.~Lou},
\newblock \bibinfo{title}{Hybrid discontinuous galerkin/finite volume method
  with subcell resolution for shocked flows},
\newblock \bibinfo{journal}{AIAA Journal}  (\bibinfo{year}{2021})
  \bibinfo{pages}{1--18}.
\bibitem[{Xiao et~al.(2005)Xiao, Honma, and Kono}]{Xiao2005}
\bibinfo{author}{F.~Xiao}, \bibinfo{author}{Y.~Honma},
  \bibinfo{author}{T.~Kono},
\newblock \bibinfo{title}{{A simple algebraic interface capturing scheme using
  hyperbolic tangent function}},
\newblock \bibinfo{journal}{International Journal for Numerical Methods in
  Fluids} \bibinfo{volume}{48} (\bibinfo{year}{2005})
  \bibinfo{pages}{1023--1040}.
\bibitem[{Xiao et~al.(2011)Xiao, Ii, and Chen}]{Xiao2011}
\bibinfo{author}{F.~Xiao}, \bibinfo{author}{S.~Ii}, \bibinfo{author}{C.~Chen},
\newblock \bibinfo{title}{{Revisit to the THINC scheme: A simple algebraic VOF
  algorithm}},
\newblock \bibinfo{journal}{Journal of Computational Physics}
  \bibinfo{volume}{230} (\bibinfo{year}{2011}) \bibinfo{pages}{7086--7092}.
\bibitem[{Xie and Xiao(2017)}]{Xie2017}
\bibinfo{author}{B.~Xie}, \bibinfo{author}{F.~Xiao},
\newblock \bibinfo{title}{{Toward efficient and accurate interface capturing on
  arbitrary hybrid unstructured grids: The THINC method with quadratic surface
  representation and Gaussian quadrature}},
\newblock \bibinfo{journal}{Journal of Computational Physics}
  \bibinfo{volume}{349} (\bibinfo{year}{2017}) \bibinfo{pages}{415--440}.
\bibitem[{Chamarthi and Frankel(2021)}]{Chamarthi2021}
\bibinfo{author}{A.~S. Chamarthi}, \bibinfo{author}{S.~H. Frankel},
\newblock \bibinfo{title}{High-order central-upwind shock capturing scheme
  using a boundary variation diminishing (bvd) algorithm},
\newblock \bibinfo{journal}{Journal of Computational Physics}
  \bibinfo{volume}{427} (\bibinfo{year}{2021}) \bibinfo{pages}{110067}.
\bibitem[{Ruan et~al.(2020)Ruan, Zhang, Tian, and He}]{Ruan2020}
\bibinfo{author}{Y.~Ruan}, \bibinfo{author}{X.~Zhang},
  \bibinfo{author}{B.~Tian}, \bibinfo{author}{Z.~He},
\newblock \bibinfo{title}{{A flux split based finite-difference two-stage
  boundary variation diminishing scheme with application to the Euler
  equations}},
\newblock \bibinfo{journal}{Computers and Fluids} \bibinfo{volume}{213}
  (\bibinfo{year}{2020}) \bibinfo{pages}{104725}.
\bibitem[{Shi et~al.(2003)Shi, Zhang, and Shu}]{Shi2003}
\bibinfo{author}{J.~Shi}, \bibinfo{author}{Y.~T. Zhang}, \bibinfo{author}{C.~W.
  Shu},
\newblock \bibinfo{title}{{Resolution of high order WENO schemes for
  complicated flow structures}},
\newblock \bibinfo{journal}{Journal of Computational Physics}
  \bibinfo{volume}{186} (\bibinfo{year}{2003}) \bibinfo{pages}{690--696}.
\bibitem[{Remacle et~al.(2003)Remacle, Flaherty, and Shephard}]{Remacle2003}
\bibinfo{author}{J.~F. Remacle}, \bibinfo{author}{J.~E. Flaherty},
  \bibinfo{author}{M.~S. Shephard},
\newblock \bibinfo{title}{{An adaptive discontinuous Galerkin technique with an
  orthogonal basis applied to compressible flow problems}},
\newblock \bibinfo{journal}{SIAM Review} \bibinfo{volume}{45}
  (\bibinfo{year}{2003}) \bibinfo{pages}{53--72}.
\bibitem[{Liska and Wendroff(2003)}]{Liska2003}
\bibinfo{author}{R.~Liska}, \bibinfo{author}{B.~Wendroff},
\newblock \bibinfo{title}{{Comparison of several difference schemes on 1D and
  2D test problems for the Euler equations}},
\newblock \bibinfo{journal}{SIAM Journal on Scientific Computing}
  \bibinfo{volume}{25} (\bibinfo{year}{2003}) \bibinfo{pages}{995--1017}.
\bibitem[{Ha et~al.(2005)Ha, Gardner, Gelb, and Shu}]{Ha2005}
\bibinfo{author}{Y.~Ha}, \bibinfo{author}{C.~L. Gardner},
  \bibinfo{author}{A.~Gelb}, \bibinfo{author}{C.~W. Shu},
\newblock \bibinfo{title}{{Numerical simulation of high mach number
  astrophysical jets with radiative cooling}},
\newblock \bibinfo{journal}{Journal of Scientific Computing}
  \bibinfo{volume}{24} (\bibinfo{year}{2005}) \bibinfo{pages}{597--612}.
\bibitem[{Zhao et~al.(2018)Zhao, Sun, Xie, and Wang}]{Zhao2018}
\bibinfo{author}{G.~Zhao}, \bibinfo{author}{M.~Sun}, \bibinfo{author}{S.~Xie},
  \bibinfo{author}{H.~Wang},
\newblock \bibinfo{title}{{Numerical dissipation control in an adaptive WCNS
  with a new smoothness indicator}},
\newblock \bibinfo{journal}{Applied Mathematics and Computation}
  \bibinfo{volume}{330} (\bibinfo{year}{2018}) \bibinfo{pages}{239--253}.
\bibitem[{Zhang et~al.(2020)Zhang, Wang, and Zhang}]{Zhang2020}
\bibinfo{author}{H.~Zhang}, \bibinfo{author}{G.~Wang},
  \bibinfo{author}{F.~Zhang},
\newblock \bibinfo{title}{{A multi-resolution weighted compact nonlinear scheme
  for hyperbolic conservation laws}},
\newblock \bibinfo{journal}{International Journal of Computational Fluid
  Dynamics} \bibinfo{volume}{34} (\bibinfo{year}{2020})
  \bibinfo{pages}{187--203}.
\bibitem[{Li et~al.(2020)Li, Wang, Zhao, Sun, Xiong, and Tang}]{Li2020}
\bibinfo{author}{L.~Li}, \bibinfo{author}{H.~B. Wang}, \bibinfo{author}{G.~Y.
  Zhao}, \bibinfo{author}{M.~B. Sun}, \bibinfo{author}{D.~P. Xiong},
  \bibinfo{author}{T.~Tang},
\newblock \bibinfo{title}{{An Efficient Low-Dissipation Hybrid Central/WENO
  Scheme for Compressible Flows}},
\newblock \bibinfo{journal}{International Journal of Computational Fluid
  Dynamics} \bibinfo{volume}{34} (\bibinfo{year}{2020})
  \bibinfo{pages}{705--730}.
\bibitem[{Peng et~al.(2021)Peng, Liu, Li, Zhang, and Shen}]{Peng2021}
\bibinfo{author}{J.~Peng}, \bibinfo{author}{S.~Liu}, \bibinfo{author}{S.~Li},
  \bibinfo{author}{K.~Zhang}, \bibinfo{author}{Y.~Shen},
\newblock \bibinfo{title}{{An efficient targeted ENO scheme with local adaptive
  dissipation for compressible flow simulation}},
\newblock \bibinfo{journal}{Journal of Computational Physics}
  \bibinfo{volume}{425} (\bibinfo{year}{2021}) \bibinfo{pages}{109902}.
\bibitem[{Li et~al.(2021)Li, Fu, and Adams}]{Li2021}
\bibinfo{author}{Y.~Li}, \bibinfo{author}{L.~Fu}, \bibinfo{author}{N.~A.
  Adams},
\newblock \bibinfo{title}{{A low-dissipation shock-capturing framework with
  flexible nonlinear dissipation control}},
\newblock \bibinfo{journal}{Journal of Computational Physics}
  \bibinfo{volume}{428} (\bibinfo{year}{2021}) \bibinfo{pages}{109960}.
\bibitem[{Don et~al.(2018)Don, Li, Wong, and Gao}]{Don2018}
\bibinfo{author}{W.~S. Don}, \bibinfo{author}{P.~Li}, \bibinfo{author}{K.~Y.
  Wong}, \bibinfo{author}{Z.~Gao},
\newblock \bibinfo{title}{{Improved symmetry property of high order weighted
  essentially non-oscillatory finite difference schemes for hyperbolic
  conservation laws}},
\newblock \bibinfo{journal}{Advances in Applied Mathematics and Mechanics}
  \bibinfo{volume}{10} (\bibinfo{year}{2018}) \bibinfo{pages}{1418--1439}.
\bibitem[{Don et~al.(2020)Don, Li, Gao, and Wang}]{Don2020}
\bibinfo{author}{W.~S. Don}, \bibinfo{author}{D.~M. Li},
  \bibinfo{author}{Z.~Gao}, \bibinfo{author}{B.~S. Wang},
\newblock \bibinfo{title}{{A Characteristic-wise Alternative WENO-Z Finite
  Difference Scheme for Solving the Compressible Multicomponent Non-reactive
  Flows in the Overestimated Quasi-conservative Form}},
\newblock \bibinfo{journal}{Journal of Scientific Computing}
  \bibinfo{volume}{82} (\bibinfo{year}{2020}) \bibinfo{pages}{1--24}.
\bibitem[{Wang et~al.(2020)Wang, Don, Garg, and Kurganov}]{Wang2020}
\bibinfo{author}{B.~S. Wang}, \bibinfo{author}{W.~S. Don},
  \bibinfo{author}{N.~K. Garg}, \bibinfo{author}{A.~Kurganov},
\newblock \bibinfo{title}{{Fifth-order A-WENO finite-difference schemes based
  on a new adaptive diffusion central numerical flux}},
\newblock \bibinfo{journal}{SIAM Journal on Scientific Computing}
  \bibinfo{volume}{42} (\bibinfo{year}{2020}) \bibinfo{pages}{A3932--A3956}.
\bibitem[{Fleischmann et~al.(2019)Fleischmann, Adami, and
  Adams}]{Fleischmann2019}
\bibinfo{author}{N.~Fleischmann}, \bibinfo{author}{S.~Adami},
  \bibinfo{author}{N.~A. Adams},
\newblock \bibinfo{title}{{Numerical symmetry-preserving techniques for
  low-dissipation shock-capturing schemes}},
\newblock \bibinfo{journal}{Computers and Fluids} \bibinfo{volume}{189}
  (\bibinfo{year}{2019}) \bibinfo{pages}{94--107}.
\bibitem[{Ii et~al.(2014)Ii, Xie, and Xiao}]{Ii2014}
\bibinfo{author}{S.~Ii}, \bibinfo{author}{B.~Xie}, \bibinfo{author}{F.~Xiao},
\newblock \bibinfo{title}{An interface capturing method with a continuous
  function: The thinc method on unstructured triangular and tetrahedral
  meshes},
\newblock \bibinfo{journal}{Journal of Computational Physics}
  \bibinfo{volume}{259} (\bibinfo{year}{2014}) \bibinfo{pages}{260--269}.
\bibitem[{Xie et~al.(2014)Xie, Ii, and Xiao}]{Xie2014}
\bibinfo{author}{B.~Xie}, \bibinfo{author}{S.~Ii}, \bibinfo{author}{F.~Xiao},
\newblock \bibinfo{title}{An efficient and accurate algebraic interface
  capturing method for unstructured grids in 2 and 3 dimensions: The thinc
  method with quadratic surface representation},
\newblock \bibinfo{journal}{International Journal for Numerical Methods in
  Fluids} \bibinfo{volume}{76} (\bibinfo{year}{2014})
  \bibinfo{pages}{1025--1042}.
\bibitem[{{Van Leer}(1979)}]{VanLeer1979}
\bibinfo{author}{B.~{Van Leer}},
\newblock \bibinfo{title}{{Towards the Ultimate Conservative Difference Scheme.
  V. A Second-Order Sequel to Godunov's Method}},
\newblock \bibinfo{journal}{Journal of Computational Physics}
  \bibinfo{volume}{32} (\bibinfo{year}{1979}) \bibinfo{pages}{101--136}.
\bibitem[{Pirozzoli(2006)}]{Pirozzoli2006}
\bibinfo{author}{S.~Pirozzoli},
\newblock \bibinfo{title}{{On the spectral properties of shock-capturing
  schemes}},
\newblock \bibinfo{journal}{Journal of Computational Physics}
  \bibinfo{volume}{219} (\bibinfo{year}{2006}) \bibinfo{pages}{489--497}.
\bibitem[{Toro et~al.(1994)Toro, Spruce, and Speares}]{Toro1994}
\bibinfo{author}{E.~F. Toro}, \bibinfo{author}{M.~Spruce},
  \bibinfo{author}{W.~Speares},
\newblock \bibinfo{title}{{Restoration of the contact surface in the
  HLL-Riemann solver}},
\newblock \bibinfo{journal}{Shock Waves} \bibinfo{volume}{4}
  (\bibinfo{year}{1994}) \bibinfo{pages}{25--34}.
\bibitem[{Toro(2009)}]{Toro2009}
\bibinfo{author}{E.~F. Toro}, \bibinfo{title}{{Riemann Solvers and Numerical
  Methods for Fluid Dynamics}}, \bibinfo{publisher}{Elsevier B.V.},
  \bibinfo{year}{2009}.
\bibitem[{Harten et~al.(1983)Harten, Lax, and van Leer}]{Harten1983a}
\bibinfo{author}{A.~Harten}, \bibinfo{author}{P.~D. Lax},
  \bibinfo{author}{B.~van Leer},
\newblock \bibinfo{title}{{On Upstream Differencing and Godunov-Type Schemes
  for Hyperbolic Conservation Laws}},
\newblock \bibinfo{journal}{SIAM Review} \bibinfo{volume}{25}
  (\bibinfo{year}{1983}) \bibinfo{pages}{35--61}.
\bibitem[{Johnsen and Colonius(2006)}]{Johnsen2006}
\bibinfo{author}{E.~Johnsen}, \bibinfo{author}{T.~Colonius},
\newblock \bibinfo{title}{{Implementation of WENO schemes in compressible
  multicomponent flow problems}},
\newblock \bibinfo{journal}{Journal of Computational Physics}
  \bibinfo{volume}{219} (\bibinfo{year}{2006}) \bibinfo{pages}{715--732}.
\bibitem[{Gottlieb(2005)}]{Gottlieb2005}
\bibinfo{author}{S.~Gottlieb},
\newblock \bibinfo{title}{{On High Order Strong Stability Preserving
  Runge-Kutta and Multi Step Time Discretizations}},
\newblock \bibinfo{journal}{Journal of Scientific Computing}
  \bibinfo{volume}{25} (\bibinfo{year}{2005}) \bibinfo{pages}{105--128}.
\bibitem[{Kurganov and Tadmor(2002)}]{Kurganov2002}
\bibinfo{author}{A.~Kurganov}, \bibinfo{author}{E.~Tadmor},
\newblock \bibinfo{title}{{Solution of two-dimensional Riemann problems for gas
  dynamics without Riemann problem solvers}},
\newblock \bibinfo{journal}{Numerical Methods for Partial Differential
  Equations} \bibinfo{volume}{18} (\bibinfo{year}{2002})
  \bibinfo{pages}{584--608}.
\bibitem[{Xu and Shu(2005)}]{Xu2005}
\bibinfo{author}{Z.~Xu}, \bibinfo{author}{C.~W. Shu},
\newblock \bibinfo{title}{{Anti-diffusive flux corrections for high order
  finite difference WENO schemes}},
\newblock \bibinfo{journal}{Journal of Computational Physics}
  \bibinfo{volume}{205} (\bibinfo{year}{2005}) \bibinfo{pages}{458--485}.
\bibitem[{Sutherland(2010)}]{Sutherland2010}
\bibinfo{author}{R.~S. Sutherland},
\newblock \bibinfo{title}{{A new computational fluid dynamics code I: Fyris
  Alpha}},
\newblock \bibinfo{journal}{Astrophysics and Space Science}
  \bibinfo{volume}{327} (\bibinfo{year}{2010}) \bibinfo{pages}{173--206}.
\bibitem[{Schneider and Robertson(2015)}]{Schneider2015}
\bibinfo{author}{E.~E. Schneider}, \bibinfo{author}{B.~E. Robertson},
\newblock \bibinfo{title}{{Cholla: A new massively parallel hydrodynamics code
  for astrophysical simulation}},
\newblock \bibinfo{journal}{Astrophysical Journal, Supplement Series}
  \bibinfo{volume}{217} (\bibinfo{year}{2015}) \bibinfo{pages}{24}.

\end{thebibliography}
\bibliographystyle{elsarticle-num-names}

\end{document}